\UseRawInputEncoding 
\documentclass[journal]{IEEEtran}
\usepackage{xcolor}
\usepackage{generic}
\usepackage[T1]{fontenc}
\usepackage{textcomp}
\usepackage[font=small,labelfont=bf]{caption}
\usepackage{amsmath,amssymb,amsfonts,mathtools}
\usepackage{amsthm}
\usepackage{algorithmicx}
\usepackage{algpseudocode}
\usepackage{algorithm}

\usepackage{graphicx}
\usepackage{textcomp}
\usepackage{framed}
\usepackage{multirow}
\usepackage{ulem}
\usepackage{soul}
\usepackage{changes}
\usepackage{amsmath}
\usepackage{mdframed}
\usepackage{hyperref}
\usepackage{epstopdf}
\usepackage{mathtools}
\usepackage{todonotes}
\usepackage{comment}
\usepackage{subcaption}
\usepackage{tabularray}
\def\BibTeX{{\rm B\kern-.05em{\sc i\kern-.025em b}\kern-.08em
    T\kern-.1667em\lower.7ex\hbox{E}\kern-.125emX}}
\usepackage[bbgreekl]{mathbbol}
\usepackage{amsfonts}
\DeclareSymbolFontAlphabet{\mathbb}{AMSb}
\DeclareSymbolFontAlphabet{\mathbbl}{bbold}

\newcommand{\calL}{{\mathcal{L}}}
\newcommand{\calU}{{\mathcal{U}}}
\newcommand{\calM}{{\mathcal{M}}}

\newcommand{\calW}{{\mathcal{W}}}

\newcommand{\sH}{\mathcal{H}}

\newcommand{\Ran}[1]{\mathcal{R}(#1)}
\newcommand{\Nul}[1]{\mathcal{N}(#1)}

\renewcommand{\AA}{\mathbb{A}}
\newcommand{\BB}{{\mathbb{B}}}
\newcommand{\HH}{{\mathbb{H}}}

\newcommand{\KK}{{\mathbb{K}}}
\newcommand{\TT}{{\mathbb{T}}}
\newcommand{\RR}{{\mathbb{R}}}
\newcommand{\MM}{{\mathbb{M}}}

\newcommand{\WW}{{\mathbb{W}}}
\newcommand{\UU}{{\mathbb{U}}}
\newcommand{\ZZ}{{\mathbb{Z}}}
\newcommand{\knl}{\mathfrak{K}}

\newcommand{\unl}{\mathfrak{u}}
\newcommand{\wnl}{\mathfrak{w}}

\newcommand{\Pwr}{\mathcal{P}}

\newtheorem{theorem}{Theorem}

\usepackage[compatibility=false]{caption}

\allowdisplaybreaks

\begin{document}

\title{Tailoring Reproducing Kernels for Optimal Control via Policy Iteration}
\author{Shengyuan  Niu$^1$, Ali Bouland$^1$, Haoran Wang$^1$, Filippos Fotiadis$^2$, \\Andrew Kurdila$^1$, Andrea L'Afflitto$^3$, Sai Tej Paruchuri$^4$, Kyriakos G.  Vamvoudakis$^5$
\thanks{$^1$S. Niu, A. Bouland, H. Wang, and A. Kurdila are with the Department of Mechanical Engineering, Virginia Tech, Blacksburg, VA, USA. Email:
{\tt\small \{syniu97, bouland, haoran9, kurdila\}@vt.edu}.}
\thanks{$^2$F. Fotiadis is with The Oden Institute for Computational Engineering \& Sciences, University of Texas at Austin, Austin, TX, 
USA. Email:
{\tt\small  ffotiadis@utexas.edu}.}
\thanks{$^3$A.  L'Afflitto is with The Grado Department of Industrial and Systems Engineering, Virginia Tech, Blacksburg, VA, USA. Email: {\tt\small a.lafflitto@vt.edu}.}
\thanks{$^{4}$S. T. Paruchuri is with the Department of Mechanical Engineering and Mechanics, Lehigh University, Bethlehem, PA, USA. Email: {\tt\small saitejp@lehigh.edu}. 
}
\thanks{$^5$ K. G. Vamvoudakis is with The Daniel Guggenheim School of Aerospace Engineering, 
Georgia Institute of Technology, Atlanta, GA, 
USA. Email:
{\tt\small kyriakos@gatech.edu}.}
\thanks{
This research was supported in part by NSF under grant Nos. CAREER CPS-$1851588$,  CPS-$2038589$, CPS-$2227185$, S\&AS-$1849198$, and SLES-$2415479$. This work was supported in part by the Office of Naval Research (ONR) under grant number  N$00014$-$24$-$1$-$2267$.}
}
\maketitle
\begin{abstract}
This paper presents a novel approach to formulating the actor-critic method for optimal control by casting policy iteration in reproducing kernel Hilbert spaces (RKHSs -- also known as native spaces).
By tailoring the reproducing kernel and  RKHS to the dynamics of the nonlinear optimal control problem,
we leverage recent advancements in characterizing error bounds from statistical and machine learning theory.
These approximations define a general strategy to select the bases of the actor-critic networks, and we formally guarantee for the first time that this basis selection procedure leads to closed-form error bounds for the individual steps of policy iteration. These bounds often have a geometric and computable form, making them potentially useful for \textit{a priori} or \textit{a posteriori} evaluation of candidate collections of scattered bases. Numerical studies subsequently provide qualitative evidence of the practical performance achieved for the full recursion using the algorithms and theory developed for the single-step error bounds. 
\end{abstract}



 
 \section{Introduction}

Autonomous systems require controllers that not only steer them to a trim point or towards a desired trajectory,  but are also efficient in terms of consumed fuel or resources.
Designing such controllers involves solving an optimal control problem, which aims to strike a balance between achieving a control objective quickly enough and minimizing the control effort.
While optimal control is relatively straightforward in simple settings characterized by linear systems and quadratic cost functions, it becomes highly challenging for general nonlinear systems.
In such scenarios, solving the Hamilton-Jacobi-Bellman (HJB) partial differential equation (PDE) is a popular theoretical approach to obtain the optimal controller.
However, due to its nonlinearity, the HJB equation is almost always impossible to solve analytically, as explained in \cite{bardi1997optimal,kamalapurkar2018reinforcement,lewis2012optimal}.

The study of optimal control theory, particularly through the formulation and solution of the HJB equation, has yielded a substantial body of literature over the past few decades. A convenient way to organize this research is to categorize approaches based on the properties of the value function, which appears in the solution of the HJB equation.

The first category introduces the concept of viscosity solutions, which were initially popularized by Crandall and Lions \cite{crandall1983viscosity}.
Accounts of this approach are provided by  \cite{bardi1997optimal,simona2024viscosity}. 
This class of methods stands out for offering a rigorous and comprehensive framework to study the existence of solutions to the HJB equation, even for solutions with value functions that exhibit relatively lower regularity.

A second category of research typically assumes higher regularity for the value function and foregoes some of the subtleties of the theory of viscosity solutions.
It focuses on developing theory and algorithms that yield practical approximately optimal controllers, often lower-dimensional and  suitable for online applications
References \cite{bertsekas2011dynamic,kamalapurkar2018reinforcement,lewis2013reinforcement} are popular examples that, for the most part,  avoid the theory of viscosity solutions and represent this latter family of approaches.
The work presented in this paper falls under this category, as does the more detailed summary of relevant research provided below.

To avoid explicitly solving a nonlinear PDE to obtain the optimal controller, the policy iteration (PI) algorithm often features prominently in this latter class. This process serves as both an explanation and a motivation for more advanced algorithms, as discussed in \cite{bertsekas2011dynamic,kamalapurkar2018reinforcement,lewis2013reinforcement}.
The policy iteration algorithm recursively evaluates and improves a given control policy until convergence is observed.  While analytically solving these PDEs remains challenging, it has become standard practice to use Galerkin approximations of the PI algorithms to design practical algorithms for synthesizing lower-dimensional, approximately optimal nonlinear controllers.

Early efforts in studying Galerkin approximations for solving linear PDEs in reinforcement learning (RL) are documented in \cite{abu2005nearly,bea1998successive,beard1997galerkin}.
In these works, the main idea is to employ a function approximation in the linear PDEs, which embodies the critic step, and subsequently train it using the Galerkin approximation or least squares methods.
Building on this foundational work, subsequent research, such as in \cite{bhasin2013novel,vamvoudakis2010online}, expanded to investigate various online implementations grounded in learning theory. Comprehensive reviews of contemporary theories that underpin many recent online and offline methods of this type can be found in \cite{lewis2013reinforcement} and \cite{kamalapurkar2018reinforcement}. 

Nevertheless, a well-known and common shortcoming of the aforementioned class is that they work effectively for lower-dimensional approximations only if the underlying function approximations have been carefully chosen.
Performance of the control strategy, therefore, relies on the assumption that the bases used are selected from some appropriate hypothesis space and that the number of basis functions is sufficient. These assumptions can be difficult to justify in practice.  
More recent studies, such as \cite{bian2021reinforcement,kalise2020robust,yang2021hamiltonian}, further emphasize the importance of investigating the impact of approximation errors on the performance of the resulting RL methods.

In light of the need for rigorous accounts of the impact of basis selection on nonlinear controller performance,
this paper leverages recent advancements in statistical and machine learning theory, or distribution-free regression.
A prominent trend in recent research on developing adaptive control schemes based on optimal control principles has been the adoption of techniques from statistical and machine learning theory.
For both discrete and continuous time dynamical systems, most of the foundational work in this field for various identification, estimation, or control problems systematically employs reproducing kernel Hilbert spaces (RKHSs) of scalar-valued functions or vector-valued RKHSs (vRKHSs) \cite[Ch. 3]{kurdila2025data}.
Notably, these efforts have emerged due to the remarkable success of RKHS theory and techniques of Gaussian processes \cite{scholkopf,williams2006gaussian} in solving numerous well-known classical learning problems for function regression based on a discrete collection of independent and identically distributed (IID) measurements. 

Among the recent efforts to synthesize optimal control theory with statistical and machine learning theory, the papers most relevant to the approach here include the very recent work in \cite{bouland2023,ehring2024hermite,ehring2025online,santin2021kernel,shengyuan2024}.
The authors' previous papers \cite{bouland2023,shengyuan2024} provided initial results on the convergence rates of the critic's approximations in policy iteration, particularly by formulating  PI in an RKHS. Both papers related the convergence rate of the critic's approximations to the fill distance of the kernel centers used for the approximation.
Reference \cite{bouland2023} focused on an offline critic formulation, while \cite{shengyuan2024} studied an online critic approach.
Although this paper only focuses on PI methods, it is more general than \cite{bouland2023} and \cite{shengyuan2024} in that it characterizes stepwise convergence rates of the PI approach for the actor and critic working together. 
 
The approach presented here is unique in comparison to \cite{bouland2023,shengyuan2024,ehring2024hermite,ehring2025online,santin2021kernel} in both methodology and generality.
We focus on deriving general approaches to tailor the definition of reproducing kernels and native spaces to obtain stepwise error bounds that closely resemble those popularized in statistical and machine learning theory for regression.  This is achieved by utilizing the feature mapping theorem to characterize the range of the bounded linear operators that define the PI algorithm. For more details on the theory underlying this approach in RKHSs, refer to \cite{saitoh2010} or \cite{saitoh}.
This overall strategy has no counterpart in \cite{ehring2024hermite,ehring2025online,santin2021kernel}, and it generalizes \cite{bouland2023,shengyuan2024}. 

By applying principles from the theory of inverse problems, we obtain stepwise error bounds for the PI algorithm that include both a sample error term and an approximation error term.
These bounds reflect the influence of both the approximation error and the contribution of numerical noise to the convergence of the PI algorithm.
In the PI algorithm, the source of the numerical noise is typically understood as the perturbations to the current iteration associated with the error introduced from approximate solutions of previous iterations.

\subsection{Notation}
 In this paper, the state space is $\mathbb{R}^n \triangleq \mathbb{R}^{n \times 1}$, and the space of control values is $\mathbb{R}^m \triangleq \mathbb{R}^{m \times 1}$. Both spaces are assumed to contain column vectors. The range of an operator $T$ is denoted $\Ran{T}$, while the nullspace is written $\Nul{T}$. 
 Let $X$ be a normed vector space and $Y$ a Banach space.
 We denote by $\calL(X,Y)$ the Banach space of bounded linear operators mapping $X\to Y$.
 We write $\calL(X)$ for $\calL(X,X)$ and
 that $X$ is continuously embedded in $Y$, and write $X\hookrightarrow Y$, whenever $X\subset Y$ and the canonical embedding $\mathcal{I}:x\in X\mapsto \mathcal{I}(x)=x\in Y$ is a bounded operator. 

This paper studies several real-valued RKHS defined over a subset $\Omega\subset \RR^n$, typically denoted in calligraphic fonts such as $\sH,\calU_i,\calW,\hat{\calW}$.
The basic properties of such spaces are discussed in Section \ref{sec:rkhs}.
The construction of these spaces always assumes that $\Omega$ is compact, and that the reproducing kernel $\knl$ that defines $\sH$  is strictly positive and bounded on the diagonal, in the sense discussed in Section \ref{sec:rkhs}.

Here, we summarize notation for $\sH$, since the corresponding notation holds for all the RKHSs in the paper.
We denote by $\mathcal{H} \triangleq \mathcal{H}(\Omega, \mathbb{R})$ the RKHS of real-valued functions over $\Omega \subset \mathbb{R}^n$ that are generated by the admissible reproducing kernel $\knl : \Omega \times \Omega \to \mathbb{R}$.
The RKHS $\mathcal{H}$ is defined as the closure of the span of the reproducing kernel functions $\knl_x : x \in \Omega \to \mathbb{R}$, where $\knl_x \triangleq \knl(\cdot, x)\in \sH$ is the reproducing kernel centered at $x \in \Omega$.  
For additional details on RKHS theory developed employing the same notation and approach to the problem as in this paper,
interested readers are referred to \cite[Ch. 3]{RKHS_MRAC_book}.

\section{Problem Formulation}
The nonlinear control problem studied in this paper is the same as in classical works such as  \cite{kamalapurkar2018reinforcement,lewis2012optimal}.
Consider a nonlinear system
\begin{align}
    \Dot{x}(t) = f(x(t)) &+ g(x(t))u(t), \quad 
    x(0) = x_0, \quad t \geq 0, \label{eqn_sys}
\end{align}
where $x(t)\in \RR^n$ is the state,
$u(t)\in \RR^m$ is the control input, $f:\RR^n\to \RR^n$, and $g:\RR^n\to \RR^{n\times m}$. It is assumed that both $f$ and $g$ are known, and we are interested in defining a feedback controller $u(t)=\mu(x(t))$ that drives any initial condition $x_0\in \Omega\subset \RR^n$ to the origin.
Among such stabilizing controllers, we seek one that is optimal with optimality defined in terms of  the  cost functional 
\begin{align}
    J(x_0, u) =& \int^{\infty}_{0} r(x(t), u(t)) \textrm{d}t  \nonumber \\ =& \int_{0}^{\infty} \left[ x^{\rm T}(t)Qx(t)+u(\tau)^{\rm T}Ru(\tau) \right] \textrm{d}t, \label{eqn_J_lqr}
\end{align}
where the symmetric matrix $Q \in \mathbb{R}^{n \times n}\succeq 0$ is the state weighting matrix, and
the symmetric matrix $R\in \mathbb{R}^{m \times m}\succ 0$ is the control effort weighting matrix. 
The value function $V_\mu$ for the admissible feedback control $\mu$ and initial condition $x_0$ is defined as 
 $ V_{\mu} ( x_0 ) = \int^{\infty}_0 r(x(\tau),\mu(x(\tau))) \textrm{d}\tau$, 
and the optimal value function $V^\star$ is defined as
   $ V^\star (x_0) = \min_{\mu\in \calU_{ad}} V_{\mu}(x_0)$.
   Here, the admissible class $\calU_{ad}$ consists of all feedback controllers $\mu:\Omega \subset \RR^n\to \RR^m$ that are stable and drive the controlled state to the origin with a finite cost \eqref{eqn_J_lqr}.
   The Hamiltonian associated with this problem is defined as 
\begin{align}
    H_{\mu}(x,\mu, V_{\mu}) &\triangleq  r(x, \mu (x))+ \nabla V^{\rm T}_{\mu}(f(x)+g(x) \mu (x)), \nonumber \\
    &\hspace{2.0em} \textrm{for all } (x,\mu, V_{\mu}) \in \RR^n \times \RR^m \times \RR,
\end{align}
and it is used to  formulate the HJB equation 
\begin{align}
    0  = \min_{\mu\in U_{ad}} H(x,\mu, V_{\mu}) = H(x^\star,\mu^\star, V_{\mu^\star}),
\end{align}
where $x^*$ is the optimal trajectory generated by the optimal control law 
\begin{align}
    \mu^\star = -\frac{1}{2} R^{-1} g^{\rm T} \nabla V_{\mu^\star },    
\end{align}
and $V^\star\triangleq V_{\mu^*}$. 

As described in detail in \cite{kamalapurkar2018reinforcement,lewis2012optimal}, the policy iteration summarized in Algorithm \ref{algorigm_PI_2}  can be used to generate a sequence of feedback controllers and value function estimates $\{(\mu^\ell,v^\ell)\}_{\ell\geq 0}$ that approximate the optimal control policy. Under suitable hypotheses,   it can be shown that  
\begin{align}
\mu^\ell \to \mu^\star ,  \quad 
v^\ell \to V^\star \triangleq V_{\mu^\star } \,  \text{ as } \ell \to \infty, \label{eqn_PI_convergence}
\end{align}
where convergence is measured in some appropriate function space norm;
a variety of choices are possible, see \cite{kundu2024policy,vamvoudakis2009online} for some examples.  
This  recursion  is defined in terms of the operators $A(\mu^{\ell-1})$ and $B$ that are defined pointwise as 
\begin{align}
    (A(\mu^{\ell-1}) v)(x) &\triangleq \left( f(x) + g(x) \mu^{\ell-1}(x) \right)^{\rm T} \nabla v(x), \label{eqn_def_operator_A} \\ 
    (Bv)(x) &\triangleq -\frac{1}{2} R^{-1} g^{\rm T}(x) \nabla v(x), \nonumber \\
    &\hspace{6.0em} \text{ for all } x\in \Omega, v\in \sH.  \label{eqn_def_operator_B}
\end{align}

\begin{algorithm}[b]
\caption{PI algorithm in an operator setting}\label{algorigm_PI_2}
    \begin{algorithmic}[1]
        \State $\mu^0(\cdot) \in \calU_{\rm ad}$,
                $\ell \gets 1$,
                $\mu^1(\cdot) \gets \text{arbitrary}\not = \mu^0(\cdot)$\Comment{Initialize the algorithm}
        \While{$l\leq \bar{\ell}$} \Comment{Repeat for all iterations, $\ell\leq \bar{\ell}$}
            
            \State Critic step: Compute $v^{\ell}$ by solving
            \begin{align}
                A(\mu^{\ell-1}) v^{\ell} &= b(\mu^{\ell-1}) \label{eqn_ideal1}
            \end{align}
            \State Actor step:
            \begin{align}
                    \mu^{\ell }(x) \gets -\frac{1}{2} R^{-1} g^{\rm T}(x) \nabla v^{\ell}(x) \label{eqn_ideal2}
            \end{align}
            \State $\ell \gets \ell+1$
        \EndWhile
    \end{algorithmic}
\end{algorithm}

\section{Characterizing the Range of Policy Iteration Operators in Native Spaces}


So far, the concept of convergence in \eqref{eqn_PI_convergence} has remained undefined. In this paper, we introduce feature mappings, feature operators, and feature spaces to define RKHSs that encompass specific combinations of derivatives present in the definition of $A\triangleq A(\mu) :\sH\to \Ran{A}\triangleq \calW$ or $B :\sH\to \Ran{B}\triangleq \calU$ in the policy iteration (PI) recursion as specified in \eqref{eqn_def_operator_A} and \eqref{eqn_def_operator_B}.  By employing these tools, we construct RKHSs that are ``well-suited'' to the optimal control problems outlined in the introduction, thereby defining the notion of convergence.

\subsection{Feature Mappings, Feature Operators, and Feature Spaces}\label{section_feature_mappings}

Let
$\HH$ be a general real Hilbert space,
which is not necessarily an RKHS.
A {feature mapping} $\Phi:\Omega\to \HH$ is a function that maps $\Omega$ into $\HH$, and
we define an associated reproducing kernel
$\knl_\Phi : \Omega\times \Omega\to \RR^+$ as
\begin{align}
    \knl_\Phi(x,y)\triangleq(\Phi(x),\Phi(y))_{\HH}, \quad \text{for all }x,y \in \Omega. \label{eqn_kernel_feature_map}
\end{align}
By definition, and as explained in Section \ref{sec:rkhs} below, we observe that $\knl_\Phi$ is an admissible real-valued kernel.
This is apparent since  the definition of the inner product $(\cdot,\cdot)_\HH$ implies  that $\knl_\Phi$ is real-valued, nonnegative, and symmetric.
Moreover, this kernel function is of positive type, as demonstrated by the following: for any $\{x_i\}_{i=1}^n\subset \Omega$ and any $\{\alpha_i\}_{i=1}^n \subset \RR$, it holds that 
\begin{align}  
    \sum_{i,j=1}^n \alpha_i \knl_\Phi(x_i,x_j)\alpha_j=&
    \left ( \sum_{i=1}^n \alpha_i\Phi(x_i), \sum_{j=1}^n \alpha_j \Phi(x_j) \right )_\HH \nonumber \\
    =&\left  \Vert 
    \sum_{i=1}^n \alpha_i \Phi(x_i)\right  \Vert ^2_\HH \geq 0.
\end{align}
Any Hilbert space $\HH$ and feature map $\Phi$ can be used to define an RKHS $\sH_\Phi$ in terms of the kernel $\knl_\Phi$.
The following theorem
describes the structural relationship between the RKHS $\sH_\Phi$ and the original Hilbert space $\HH$. 

\begin{theorem}[{\cite[Lemma 5]{xu2007refinable}, \cite[Prop. 1]{carmeli2005reproducing}}]\label{theorem_feature}
    Let $\HH$ be a real Hilbert space,
    $\Phi : \Omega\to \HH$ be a feature map,
    $\knl_\Phi$ be the associated kernel defined as in \eqref{eqn_kernel_feature_map}, and
    $\sH_\Phi$ be the real RKHS induced by the kernel $\knl_\Phi$. 
    Then,
    \begin{align}
        \sH_\Phi = &\left \{
    f: \Omega \to \RR  :  f=Fh   = (h,\Phi(\cdot))_\HH \text{ for some } h\in \HH
    \right \}\nonumber \\ \triangleq& \Ran{F}, 
    \end{align}
    with  
    \begin{align}
     &\Vert f \Vert _{\sH_\Phi} \nonumber \\
     &\, \triangleq \inf \left \{
     \Vert h \Vert _\HH  :  f\triangleq Fh\triangleq (h,\Phi(\cdot))_\HH \text{ for some } h\in \HH
    \right \}.
    \end{align}
    The feature operator $F:\HH\to \sH_\Phi$ defined so that
    $Fh\triangleq (\Phi(\cdot),h)_\HH$ for any $h\in \HH$ is a partial isometry  from $\HH$ onto $\sH_\Phi$.
    The operator $F^*F$ is the $\HH$-orthogonal projection onto $\Nul{F}^\perp$.  
\end{theorem}
 
Since $F:\HH\to \sH_\Phi$ is a partial isometry
with initial space $\Nul{F}^\perp$ and final space $\Ran{F}=\sH_\Phi$,
it follows from Theorem 4.34 of \cite{weidmann} that
\begin{align*}
    F^*F\triangleq \Pi_{\Nul{F}^\perp}\quad \text{ and } \quad 
    FF^*\triangleq I_{\sH_\Phi}. 
\end{align*}
Here, $\Pi_{\Nul{F}^\perp}$ denotes the $\HH$-orthogonal projection onto $\Nul{F}^\perp$. 
If we choose, for example, the feature operator  $F$ to be $A$, $A^*$, or $B$,  
then these identities can be crucial to studying the properties of Galerkin approximations of the associated operator equations.

\subsection[A Study of the Range of A, A*, and A*A]{A Study of the Range of $A\triangleq A(\mu),A^*\triangleq A^*(\mu)$ and $A^*A\triangleq A^*A(\mu)$}
 In this section, we study the operators
 $A\triangleq A(\mu),A^*\triangleq A^*(\mu)$ and $A^*A\triangleq A^*A(\mu)$ for a fixed feedback control law $\mu$.
 The operator $A$ is defined pointwise with 
\begin{align}
    \left ( Av\right)(x)&=\left ( A(\mu)v\right)(x)\nonumber\\
    &= \sum_{k=1}^n\left( f_k(x)+\sum_{j=1}^m  g_{kj}\mu_j(x) \right ) D^{e_k}v(x) \nonumber\\
    &= \sum_{k=1}^n a_k(x)D^{e_k}v(x),
\end{align}
where
$
    a_k(x) \triangleq f_k(x) + \sum_{j=1}^m g_{kj}\mu_j(x)$, for all $ k \in \{1, \ldots, n\} \text{ and } x \in \Omega$.

\subsubsection{A study of the case when \texorpdfstring{$A\in \calL(\sH,\calW), \text{ with } \calW=\Ran{A}$}{A in L(H,W) with W=Ran(A)}}\label{section_approx_A_H_W}
We assume that $\Omega\subset \RR^n$ is compact, and that the reproducing kernel $\knl:\Omega \times \Omega \to \RR$  that defines $\sH$ satisfies the regularity condition $\knl\in C^{2s}(\Omega\times \Omega)$ for some integer smoothness index $s\geq 1$. It follows from Theorem \ref{theorem_zhou} that
\begin{align}
    (D^{e_i}v)(x)=\left (D^{(e_i,0)}\knl(x,\cdot),v \right)_\sH \quad \text{ for every } v \in \sH,
\end{align}
where the derivative $D^{e_i}v$ is defined as described in Section \ref{sec:rkhs}. 
We choose the feature mapping $\Phi:x\mapsto \Phi(x)\in \sH$, and define the associated feature operator $F$ to be
\begin{align}
    &(Fv)(x) \triangleq (\Phi(x),v)_\sH =(Av)(x) \nonumber \\
    &\quad = \left (\sum_{k=1}^n (f_k(x) + \sum_{j=1}^m g_{kj}(x)\mu_j(x))D^{(e_k,0)}\knl(x,\cdot),v \right )_{\sH}, \label{eqn_PHI1}
\end{align}
for all $x\in \Omega$ and $v\in \sH$. 
The fact that $\Phi(x)\in \sH$ is a consequence of Theorem \ref{theorem_zhou}:
for each $x \in \Omega$,
$\Phi(x)$ is a linear combination of terms in the form $D^{(e_k,0)}\knl(x,\cdot)$, each of which is contained in $\sH$.
With this choice of the feature mapping $\Phi(x)$ and feature operator $F$,
it holds that $A\equiv F:\sH\to \calW\triangleq \Ran{F}\equiv \Ran{A}$.
The space $\calW\triangleq \Ran{A}$ is an RKHS with underlying reproducing kernel $\wnl:\Omega\times \Omega \to \RR$  guaranteed by the feature mapping Theorem \ref{theorem_feature}, which then takes the form  
\begin{align}
    \wnl(x,y)=&\left (\Phi(x),\Phi(y) \right)_\sH \nonumber \\
    =&\left (\sum_{k=1}^n a_k(x) D^{(e_k,0)}\knl(x,\cdot),\right.  \nonumber \\
    & \quad \quad \quad \quad \quad \quad \left. \sum_{\ell=1}^n a_\ell(y) D^{(e_\ell,0)}\knl(y,\cdot)  \right )_\sH  \nonumber \\
    =&\left (f(x)+g(x)\mu(x) \right)^{\rm T} \nonumber\\& \quad \quad \quad [D^{(e_k,e_\ell)}\knl(x,y)] (f(y)+g(y)\mu(y)), \label{eqn_wnldefinition}
\end{align}
for all $x,y\in \Omega$. 

The adjoint $A^*:\calW \to \sH$, interpreting $A$ as the operator $A:\sH\to \calW$  in this section, can be computed noting that,
for all $v\in \sH,w\in \calW$,
{\small 
\begin{align}
    \left (Av,w \right)_\calW &= \left (\sum_{k=1}^n {a_k(\cdot)D^{e_k}v(\cdot),w(\cdot)} \right )_{\calW} \nonumber \\
    &= \left (\sum_{k=1}^n a_k(\cdot)\left ( D^{(e_k,0)}\knl(\cdot,\circ),v(\circ)\right)_\sH,w(\cdot) \right )_{\calW}, \nonumber \\
    &=\left (v(\circ), \left (\sum_{k=1}^n a_k(\cdot) D^{(e_k,0)}\knl(\cdot,\circ),w(\cdot)\right)_\calW \right )_{\sH} \nonumber \\
    &=\left (v,A^*w\right)_\sH.  \label{eqn_temp_01}
\end{align}
}
In the equation above we use  $(\cdot)$ to denote the spatial variable for the calculation of inner products in $\calW$, and we use $(\circ)$ to be a placeholder for the spatial variable for inner products in $\sH$. 


\subsubsection{A study of the case when  \texorpdfstring{$A\in \calL(\sH,L^2(\Omega))$}{A ∈ L(H, L2(Ω))}, with \texorpdfstring{$\Ran{A}\hookrightarrow L^2(\Omega) $}{R(A) → L2(Ω)}}
\label{sec:A_H_to_L2}
The previous section introduced a feature operator $F:\sH \rightarrow \calW=\Ran{F}$, with $\calW$ being a suitable RKHS. 
%
By virtue of Theorem \ref{theorem_zhou},  which guarantees the continuous embedding of $\sH\hookrightarrow C^s(\Omega)$, we know in particular that $Av\in \calW \hookrightarrow C(\Omega)\hookrightarrow L^2(\Omega)$.
Thus, in contrast to the last section, in this section,
we interpret $A$ as a mapping from $\sH$ to $L^2(\Omega)$, so that  $A:\sH\to L^2(\Omega)$. In this interpretation, we have 
\begin{align}
    A: \sH \to \calW \overset{\mathcal{I}}{\hookrightarrow} C(\Omega) \overset{\mathcal{J}}{\hookrightarrow} L^2(\Omega),
\end{align}
with $\mathcal{I}$ and $\mathcal{J}$ being canonical embeddings. 
The operators $\mathcal{J}$, $\mathcal{I}$, and $A:\sH\to \calW$ are bounded and linear.
Since $\Omega$ is compact, it is well-known that the Arzela-Ascoli theorem implies that  $\mathcal{J}$ is compact. Hence, the operator $A:\sH\to L^2(\Omega)$, being a product of compact and bounded operators, is itself compact.

We calculate the adjoint $A^*:L^2(\Omega)\to \sH$ of the operator $A:\sH\to L^2(\Omega)$ as  
\begin{align}
    \left (Av,q\right)_{L^2(\Omega)} &=\left ( \sum_{k=1}^n a_k D^{e_k}v,q\right)_{L^2(\Omega)} \nonumber \\
    &= \int_\Omega \sum_{k=1}^na_k(\xi) \left(D^{e_k}\knl(\xi,\cdot), v  \right)_\sH q(\xi)\textrm{d}\xi \nonumber \\
    &= \left ( v, \sum_{k=1}^n\int_\Omega  a_k(\xi)D^{(e_k,0)}\knl(\xi,\cdot) q(\xi)\textrm{d}\xi\right)_{\sH} \nonumber\\ &= (v,A^*q)_{\sH}.
\end{align}
Thus, viewing the operator $A\in \calL(\sH,L^2(\Omega))$, we have that
\begin{align}
    (A^*q)(x)&\triangleq 
\int_\Omega {\sum_{k=1}^n a_k(\xi)D^{(e_k,0)}\knl(\xi,x) }q(\xi)\textrm{d}x \nonumber \\
&\triangleq (\hat{\Phi}(x),q)_{L^2(\Omega)}
\end{align}
for each $q\in L^2(\Omega)$. 

Similarly to the last section, we now define the feature mapping $x\mapsto  \hat{\Phi}(x)\in L^2(\Omega)$ and the associated feature operator $\hat{F}=A^*:L^2(\Omega)\to \sH$ as 
\begin{align}
    (\hat{F}q)(x)=(A^*q)(x)&=\left ( \hat{\Phi}(x),q\right)_{L^2(\Omega)}.
\end{align}
We define $\hat{\calW}\triangleq \Ran{\hat{F}}\triangleq \Ran{A^*}$, which is an RKHS with the kernel 
\begin{align}
    \hat{w}(x,y)&\triangleq \left (\hat{\Phi}(x),\hat{\Phi}(y) \right)_{L^2(\Omega)} \nonumber \\
    &=\sum_{k,\ell=1}^n \int_\Omega a_k(\xi)D^{(e_k,0)}\knl(\xi,x) a_\ell(\xi)D^{(e_l,0)}\knl(\xi,y)\textrm{d}\xi\notag \\
    &=\int_\Omega (f(\xi)+g(\xi)\mu(\xi))^{\rm T} \nonumber \\
    & \quad \cdot \nabla_{\xi}\knl_x(\xi)(\nabla_\xi\knl_y(\xi))^{\rm T}(f(\xi)+g(\xi)\mu(\xi)) \textrm{d}\xi. \label{eqn_kernelC3}
\end{align}
This is the RKHS and associated kernel studied in \cite{bouland2023}, where we consider the solution of the operator equation 
$
A^*A \hat{v}=A^*b
$  
as an approximation $\hat{v}=A^\dagger b=(A^*A)^\dagger A^*b$ of the solution to the equation $Av=b$.  

%
%
\subsection{The Range of  \texorpdfstring{$B$}{B}}\label{section_feature_Uspace}

By definition, it holds that
$\mu=Bv=-\frac{1}{2}R^{-1}g^{\rm T} \nabla v$ in the PI iteration 
for $v\in \sH$.
However, similarly to the procedure in Section \ref{section_approx_A_H_W},
we apply Theorem \ref{theorem_zhou} to define a feature mapping to characterize the range of $B$.
To begin, we let $\mu\triangleq \{\mu_1,\ldots,\mu_m\}^{\rm{T}}$, $Bv\triangleq [B_1v,\ldots,B_mv]^{\rm{T}}\in \calU=\bigtimes_{i=1}^m \mathcal{U}_i$, so that each $B_i:\sH\to \mathcal{U}_i$.   We construct a feature mapping $\Phi_i(x)\in \sH$ using the sequence of identities 
\begin{align}
    \mu_i(x)&\triangleq (B_iv)(x) \nonumber\\&=-\frac{1}{2}\sum_{j=1}^m\sum_{k=1}^n R^{-1}_{ij}g_{kj}(x)\underbrace{D^{(e_k,0)}v(x)}_{\left ( D^{(e_k,0)}\knl_x(\cdot), v\right)_{\sH}} \nonumber \\ 
    &=\left (-\frac{1}{2} \sum_{j=1}^m\sum_{k=1}^n R^{-1}_{ij}g_{kj}(x) D^{(e_k,0)}\knl(x,\cdot),v\right)_\sH \nonumber\\
    &=\left (\Phi_i(x), v \right)_\sH,
\end{align}
where $R_{ij}^{-1}\triangleq (R^{-1})_{ij}$. 
We define the associated feature operator $F_i:\sH\to \calU_i$ from the identity $
    (F_i v)(x) = \left (\Phi_i(x),v \right )_\sH$, 
so that $F_iv=B_iv=\left ( \Phi_i(\cdot),v\right )_{\sH}$. We then define the reproducing kernel $\unl_i:\Omega \times \Omega \to \RR$  that defines $\calU_i$ to be 
\begin{align}
    \unl_i(x,y)&\triangleq \left (\Phi_i(x),\Phi_i(y) \right)_\sH, \nonumber\\
    &=\frac{1}{4} \left ( \sum_{j=1}^m\sum_{k=1}^n R^{-1}_{ij}g_{kj}(x) D^{(e_k,0)}\knl(x,\cdot),\right. \nonumber\\ & \quad \quad \quad \quad \quad \left. \sum_{s=1}^m\sum_{t=1}^n R^{-1}_{is}g_{ts}(x) D^{(e_t,0)}\knl(y,\cdot)\right)_\sH,\notag\\
    &=\frac{1}{4}(R^{-1}_{i})^{\rm T} g^{\rm T}(x)\left [D^{(e_k,e_t)}\knl(x,y) \right ]g(y)R^{-1}_{i}, \label{eqn_defn_unl_i} \\
    &\quad \quad \text{ for all } x,y\in \Omega, \nonumber
\end{align}
where $R^{-1}_i\triangleq \left[R^{-1}_{i,1},\ldots,R^{-1}_{i,m}\right]^{\rm T}\in \RR^m$.

\section{Galerkin Approximations for Policy Iteration}
\subsection{Approximations of the Critic Update Step}
This section discusses how the tailored RKHSs 
$\calW\triangleq \Ran{A(\mu)}$, 
$\calU\triangleq \Ran{B}$, or $\hat{\calW}\triangleq \Ran{A^*}$ can be used to define Galerkin approximations for different choices of finite-dimensional subspaces in practical realizations of the PI algorithm. 
Here, we make use of some properties of Galerkin approximations as they appear in the theory of inverse problems \cite{kress1989linear,kirsch2011introduction}, which may be less familiar to the controls community.
For this reason, the general form and analysis of Galerkin approximations are, for convenience, summarized in Section \ref{section_galerkinapproximations} below. 

Since the calculations during the critic update hold for a fixed choice feedback control $\mu^{\ell-1}$,
we simply write $A\equiv A(\mu^{\ell-1})$ in this section.
Throughout this section on the critic,
we also simply write $N\triangleq N(\ell)$, suppressing the notation that emphasizes that, in general, the dimension can depend on the iteration number $\ell$.
We choose a Banach space $Y$, referred to as the universe, that satisfies 
$ H_N\subseteq \sH\to  \Ran{A}\triangleq \calW\hookrightarrow Y$.
Admitting the design or choice of $Y$ enables more flexibility in varying uncertainty classes in inverse problems; see references \cite{kirsch2011introduction,kress1989linear}.  
In our specific algorithms below, the choice of $Y$ is either $\calW$ or $L^2(\Omega)$.  Each will yield its own specific algorithm, with its own specific performance bounds. 

For these choices of $Y$, it is assumed that we also define a family of bounded projection operators $P_N: Y \to Y_N$, with  $P_N^2=P_N$ and $P_N\in \calL(Y)$.  
With these choices, the ideal Petrov-Galerkin approximation $v_N\triangleq v_{N(\ell)}\in H_N\triangleq H_{N(\ell)}$    of the solution $v\triangleq v^\ell$ of the noise-free operator equation, which appears in \eqref{eqn_ideal1} and is repeated for convenience here, 
\begin{align}
Av=b,\label{eqn_idealoperator}
\end{align}
satisfies $
    v_{N}={\left (P_NA|_{H_N} \right )^{-1}P_N}b\triangleq {G_N}b$.  
Analogously, we denote by $v_N^\delta\triangleq G_Nb^\delta$ the Galerkin approximation of the solution $v^\delta\triangleq v^{\ell,\delta}$ to the ideal noisy operator equation 
\begin{align}
    Av^\delta = b^\delta. \label{eqn_idealnoisyoperator}
\end{align}
We interpret the noisy data $b^\delta$ as representing the perturbation to the current iteration of errors arising from approximate solutions in previous steps. 

The coordinate representation of \eqref{eqn_idealnoisyoperator} can be obtained using the unique representation of $P_N$ as  
$P_Ny=\sum_{k=1}^N \lambda_{N,k}(y) y_k$,
where $\lambda_{N,k}:Y\to \RR$ is a bounded linear coordinate functional associated with the basis $\{y_1,\ldots,y_N\}\in Y$ of $Y_N$.  
Using this representation of $P_N$, we obtain the  coordinate representation  of the operator equation in the critic step as 
$
\AA_N \Theta_N = \BB_N, 
$ 
with $\AA_N\in \RR^{N\times N}$ given by $\AA_{kj}=\lambda_{N,k}(A\knl_j)$ and  $\BB_N=\left[\lambda_{N,1}(b),\ldots,\lambda_{N,N}(b)\right]\in \RR^N$.
We consider two specific choices of finite-dimensional subspaces and universe  $Y$ in the next two sections. 
\subsubsection[Critic Method  (C1), when Y=W=Ran(A)]{Critic Method   (C1),  when $Y=\calW=\Ran{A}$}  
\label{section_critic_C1}
We choose $P_N$ to be the $\calW$-orthogonal projection $\Pi_{W_N}$ of $\calW$ onto $W_N$, and we select
    \begin{align}
        H_N&\triangleq \text{span}\{\knl_{\xi_1},\ldots,\knl_{\xi_N} \} \subset \Nul{A}^\perp \subset \sH, \label{eqn_critic_C1_HN}\\
        W_N&\triangleq Y_N\triangleq \text{span}\{\wnl_{\xi_1},\ldots,\wnl_{\xi_N} \}, \label{eqn_critic_C1_WN}
    \end{align}
    where $\knl_{\xi_\ell}\triangleq \knl(\cdot,\xi_\ell)$ and $\wnl_{\xi_\ell}\triangleq \wnl(\cdot,\xi_\ell)$ for centers $\xi_\ell\in \Xi_N\subset \Omega$, and  the reproducing kernels $\knl$ and $\wnl$ in \eqref{eqn_wnldefinition} define $\sH$ and $\calW$, respectively.  
    In this case, the Petrov-Galerkin equation is equivalent to the variational equations where we seek $v_N\in H_N$ such that 
    $
    \left (Av_N-b,\wnl_k\right)_\calW =0 \quad \text{ for all } 1\leq k \leq N$. 
    The coordinate realizations are written in terms of the matrices  
    \begin{align}
    \AA_{kj}&=(A\knl_{\xi_j})(\xi_k) =(f(\xi_k)+g(\xi_k)\mu(\xi_k))\nabla \knl_{\xi_j}(\xi_k),\\
    \BB_{k}&=b(\xi_k).
    \end{align}
    We refer to these choices as defining \textit{collocation equations} in the pair  $(\sH,\calW)$.

\subsubsection{Critic Method  (C2), when \texorpdfstring{$Y=L^2(\Omega)$}{Y=L2(Omega)}}\label{section_critic_C2}
Initial studies of this approach can be found in \cite{bouland2023}. In this case, we specifically choose the basis  and finite-dimensional subspace 
    \begin{align}
        Y_N&\triangleq \text{span} \left\{y_j\triangleq A\knl_{\xi_j}  :  \knl_{\xi_j}\in H_N \right \} \subset Y\triangleq L^2(\Omega), 
    \end{align}
    and the operator $P_N$ is selected to be the $L^2(\Omega)$-orthogonal projection onto $Y_N$. The governing equations for the Galerkin approximation $v_N$ take the variational form $
    \left (Av_N-b,A{\knl_{\xi_k}}   , \right)_{L^2(\Omega)} =0 \quad \text{ for all } 1\leq k \leq N$. 
    This setting can also be interpreted   as a Bubnov-Galerkin approximation $v_N$ of the solution of the operator equation $
    A^*Av=A^*b\in \sH$, 
    so here  $v_N$ can alternatively interpreted as an approximation of the pseudoinverse $v_N\approx  A^\dagger b=(A^*A)^\dagger A^* b$.  In this interpretation the range of $A^*A$ is contained in the RKHS $\hat{\calW}$ defined by the reproducing kernel $\hat{\wnl}$ in \eqref{eqn_kernel_feature_map}. 
    In this case, the  coordinate  realizations become 
    \begin{align}
    \AA_{kj}&=\int_\Omega (A\knl_{\xi_j})(x)(A\knl_{\xi_k})(x)\textrm{d}x, \nonumber \\
    &=\int_\Omega \nabla \knl_{\xi_k}(x)^{\rm T}(f(x)+g(x)\mu(x))^{\rm T}\nonumber \\
    & \quad \quad \quad (f(x)+g(x)\mu(x))\nabla \knl_{\xi_j}(x)\textrm{d}x,  \label{eqn_AA_C2}\\
    \BB_{k}&=\int_\Omega \knl(\xi_k,x)b(x)\textrm{d}x. \label{eqn_BB_C2}
    \end{align}
     Note that, in contrast to the Critic Method (C1) above, practical calculations based on these realizations require their further approximation using quadrature; see \cite{bouland2023}.
     Implications of this observation are discussed in the numerical examples. 

In the following sections, we discuss how explicit error bounds can be derived for each of these variants of Galerkin approximations for the critic step of the PI iteration.
   
\subsection{Error Bounds for the Critic}
\subsubsection{Error Bounds for the  Critic Method  (C1)}
\label{section_error_critic_C1}
We start with Method $(C1)$ since the analysis of the single step error in the  PI method is somewhat simpler than that for Method  $(C2)$. In general, we have 
 $\Vert G_N \Vert \to \infty $,
for the Galerkin approximation operator when $A$ is a compact operator \cite{kress1989linear}.  However, since we have identified  $A$ in this case as a feature operator from $\sH$ onto  $\calW$,  the operator $A:\sH\to \calW$ is a partial isometry. The restriction of $A$ to $\Nul{A}^\perp$ is an isometry, and it has a bounded inverse  $(A|_{\Nul{A}^\perp})^{-1}:\calW\triangleq \Ran{A}\to \sH$.  It consequently is not a compact operator, unless $\Nul{A}^\perp$ is finite-dimensional.
As we will see in the proof of the next theorem, this implies that we have the uniform bound $ \Vert G_N \Vert  \leq C_1$ for some positive constant $C_1>0$ in the problem at hand. 

We need one more preparatory result to derive error bounds for Method  (C1). It is based on the proof of Theorem 11.23 of \cite{wendland}. We define the integral operator $L_{\knl}:L^2(\Omega)\to \sH$ in terms of the kernel $\knl(x,y)$ for the space $\sH$ via the expression
\begin{align}
    (L_\knl p)(x)\triangleq \int_\Omega \knl(x,y)p(y)\textrm{d}y    
\end{align}
for $p\in L^2(\Omega)$. 
We say that $f\in \sH$ satisfies a regularity assumption induced by $L_\knl$ whenever $f=L_\knl p$ for some function $p\in L^2(\Omega)$. The important result that we use in the proof below is that 
$(g,L_\knl p)_{\sH}=(g,p)_{L^2(\Omega)}$ for all $g\in \sH, p\in L^2(\Omega)$. This is easy to show when we interpret the operator $L_\knl$ in terms of a Bochner integral. This holds since then  for any $g\in \sH$ and $p\in L^2(\Omega)$ we have
{\begin{align}
(L_\knl p,g)_\sH &= \left ( \int_\Omega \knl_y(\cdot)p(y)\textrm{d}y,g \right)_\sH \nonumber \\
&= \int_\Omega(\knl_y,g)_\sH p(y)\textrm{d}y \nonumber \\
&=\int_\Omega g(y)p(y)\textrm{d}y=(p,g)_{L^2(\Omega)}. 
\end{align}
}
Using this uniform bound on $G_N$ and the above regularity assumption, we have the following error bounds that describe a single iteration of the critic step.

\begin{theorem}\label{theorem_method C1}
Choose the finite dimensional subspaces $H_N\subset \Nul{A}^\perp \subset  \sH$ and $W_N\subset \calW$ as in \eqref{eqn_critic_C1_HN} and \eqref{eqn_critic_C1_WN},  and further suppose  $v^\delta_N\in \sH$   is   the Galerkin approximation obtained using  Method  (C1) for the    ideal noisy  operator \eqref{eqn_idealnoisyoperator}. 
Further, assume that the solution  $v$ of the ideal operator \eqref{eqn_idealoperator} without noise 
satisfies the regularity assumption $v=L_\knl p$ for some $p\in L^2(\Omega)$. 
For a convergent Galerkin approximation,  there exist constants $C_1,C_2>0$ such the error bound
\begin{align}
     \Vert v&-v_N^\delta \Vert _{L^\infty(\Omega)}\leq  \Vert v-v_N^\delta \Vert _{\sH}, \nonumber \\
     &\leq C_1 \delta + C_2 \Vert L_\knl^{-1}b \Vert _{L^2(\Omega)} \nonumber \\
     &\quad \cdot \sup_{\xi\in \Omega} \sqrt{\knl(\xi,\xi)-\knl_{\Xi_N}^\mathrm{T}(x)\KK_N^{-1}\knl_{\Xi_N}(x)} \nonumber \\
    &=C_1 \delta + C_2\sup_{\xi\in \Omega} \mathcal{P}_{H_N}(\xi)  \Vert L^{-1}_\knl v \Vert _{L^2(\Omega)}
\end{align}
is verified
with $\delta\triangleq  \Vert b-b^\delta \Vert _\calW$,
$\knl_{\Xi_N}(x) = \left[\knl_{\xi_1}(x),\ldots,\knl_{\xi_N}(x)\right]^\mathrm{T}\in \RR^N$,  $\KK_N\triangleq [\knl(\xi_i,\xi_j)]\in \RR^{N\times N}$, and $\mathcal{P}_{H_N}$ the power function (see Section \ref{sec:rkhs}) of the subspace $H_N$ in the RKHS $\sH$. 
\end{theorem}
\begin{proof}
    The overall form of the bound in this theorem follows by applying the analysis and symbols summarized in Section \ref{section_galerkinapproximations} below, as well as  
    \cite[Th. 17.6]{kress1989linear},
    with $T\triangleq A$, $X\triangleq \Nul{A}^\perp$, $  P_N\triangleq \Pi_{W_N}$, 
        $Y=Z\triangleq \calW=\Ran{A}$, $X_N\triangleq H_N\subset \Nul{A}^\perp$, and 
        $Y_N=Z_N\triangleq W_N$.  
    Here, $\Pi_{W_N}$ denotes the $\calW$-orthogonal projection onto $W_N\subset \calW$.  The bound in this theorem is based on refining the baseline, well-known estimate 
    \begin{equation}
 \Vert v-v_N^\delta \Vert _\sH \leq  \Vert G_N \Vert \delta + C \inf_{g_N\in H_N} \Vert v-g_N \Vert _{\sH}, \label{eqn_baseline_estimate}  
\end{equation} 
which, again, is  summarized in Section \ref{section_galerkinapproximations} below. 
    It is easy to check that using Theorem 17.6  of \cite{kress1989linear} in this case there is a constant $C_1>0$ that yields $ \Vert G_N \Vert \leq C_1$. This is not surprising in view of the fact that $A|_{\Nul{A}^\perp}$ is an isometry in the case at hand. With these choices, we have
 $
 G_N\triangleq \left ( \Pi_{W_N}A|_{H_N}\right)^{-1}\Pi_{W_N}:\calW \to H_N\subset \sH.  
 $ 
 Theorem 17.6  of \cite{kress1989linear} yields 
 $
  \Vert G_N \Vert _{\calL(\calW,\sH)} \leq C  \Vert \Pi_{W_N} \Vert  \sup_{z_N\in W_N, \Vert z_N \Vert _\calW=1} \Vert z_N \Vert _{\calW}\triangleq C_1.
 $
    The rightmost  term in the stepwise bound above 
  follows from the identity
    \begin{align}
         \Vert (I-\Pi_{H_N} )v \Vert _{\sH} &\leq  \Vert \mathcal{P}_{H_N} \Vert _{L^2(\Omega)} \Vert p \Vert _{L^2(\Omega)}, 
    \end{align}
which is derived from an intermediate step in the proof of \cite{wendland}. Here, for completeness,  we reproduce the lines from this proof that are relevant.  We have 
\begin{align}
     \Vert (I-\Pi_{H_N})v \Vert ^2_\sH&=((I-\Pi_{H_N})v,(I-\Pi_{H_N})v)_{\sH} \nonumber\\
    &= \left((I-\Pi_{H_N})v,v\right)_{\sH} \nonumber \\
   &=((I-\Pi_{H_N})v,L_\knl p)_{\sH} \nonumber\\
   &= ((I-\Pi_{H_N})v,p)_{L^2(\Omega)}, \nonumber \\
  &\leq  \Vert (I-\Pi_{H_N})v \Vert _{L^2(\Omega)} \Vert p \Vert _{L^2(\Omega)}.\nonumber
\end{align}
However, as discussed in Section \ref{sec:rkhs} below,
the pointwise error bound 
\begin{align}
    |E_x(I-\Pi_{H_N})v|\leq \mathcal{P}_{H_N}(x) \Vert (I-\Pi_{H_N})v \Vert _\sH
\end{align}
is always verified, which yields 
\begin{align*}
 \Vert (I-\Pi_{H_N})v \Vert ^2_{L^2(\Omega)} &\leq \int_\Omega |\mathcal{P}_{H_N}(\xi)|^2  \Vert (I-\Pi_{H_N}) v \Vert _\sH^2 \textrm{d}\xi\\
 &= \Vert \mathcal{P}_{H_N} \Vert _{L^2(\Omega)}^2  \Vert (I-\Pi_{H_N})v \Vert _\sH^2.
\end{align*}
We obtain the final form of the bound  in the theorem using the relationship  
\begin{align}
 \Vert \mathcal{P}_{H_N} \Vert ^2_{L^2(\Omega)}&=\int_\Omega (\mathcal{P}_{H_N}(\xi))^2 \textrm{d}\xi\leq |\Omega|\sup_{\xi\in \Omega}\left|\mathcal{P}_{H_N}(\xi) \right|^2,  
\end{align}
with $C_2\triangleq \sqrt{|\Omega|}$ and $|\Omega|$ is the Lebesgue measure of $\Omega\subset \RR^n$. 
\end{proof}

\subsubsection{Error Bounds for Critic Method  (C2)}\label{section_criticC3}
Now, we turn to the analysis of error bounds for the critic Method  (C2).
The analysis of error bounds for this case is studied in \cite{bouland2023} for the noise-free setting.
In this section, we extend the analysis to the situation when the data is subject to perturbation, which is due to the effect of approximate solutions of previous steps.    

Here, in the language and symbols used in the appendix,
we  choose $T\triangleq A$, $X\triangleq \Nul{A}^\perp$, $Z\triangleq \calW=\Ran{A}$, $Y\triangleq L^2(\Omega)$, 
        $X_N\triangleq H_N\subset \Nul{A}^\perp$, $Z_N\triangleq (W_N,\calW)$, $Y_N\triangleq (W_N,L^2(\Omega))$, and $
        P_N\triangleq \Pi_{Y_N}\triangleq \Pi^{L^2(\Omega)}_{Y_N}$. 
The matrix realizations described above for Method  $(C2)$ are equivalent to selecting the Galerkin operator as $
G_N\triangleq \left (\Pi_{Y_N}A|_{H_N} \right)^{-1}\Pi_{Y_N}$, 
where $\Pi_{Y_N}$ is the $L^2(\Omega)$-orthogonal projection onto $(Y_N,L^2(\Omega))$.
We newly leverage the basic estimate in \eqref{eqn_baseline_estimate}, but
now we must use more care in bounding $ \Vert G_N \Vert $.
The operator $A:\sH\to L^2(\Omega)$ is compact in this case, as discussed in \cite{bouland2023} and summarized in Section \ref{sec:A_H_to_L2}.
We expect that $ \Vert G_N \Vert \to \infty$ as $N\to \infty$.  Indeed, as discussed in Section \ref{section_galerkinapproximations}, we use an inverse inequality for the pair of spaces $\calW\hookrightarrow L^2(\Omega)$. 

\bigskip

\paragraph{\underline{Single Step Error Bounds Using  Bernstein Inequalities}}
In this case, we say that the pair $\calW\hookrightarrow L^2(\Omega)$ and bases $W_N$ satisfy an inverse, or Bernstein,  inequality \cite{devore1993constructive} with rate function $R_I(N)\to 0$ if there is a constant $C_I>0$ such that 
\begin{equation}
 \Vert w_N \Vert _{\calW} \leq C_I\frac{1}{R_I(N)}  \Vert w_N \Vert _{L^2(\Omega)} \quad \text{ for all } w_N\in W_N. \label{eqn_BIneqC3}
\end{equation}
We can argue as  in Theorem 17.2 of \cite{kress1989linear} that,  when  this inverse inequality holds, we have 
 $\Vert G_N \Vert \leq C_I \frac{1}{R_I(N)}. $
This expression provides a bound that is needed to refine the basic estimate above.


\begin{theorem}
\label{theorem_method_C2}
Let $v^\delta_N$   be  the Galerkin approximation obtained using  Method  ($C2$) for the noisy   ideal operator equation 
$
Av^\delta=b^\delta \in L^2(\Omega)$.  
Further assume that the solution $v$ of the noise-free equation \eqref{eqn_ideal1} satisfies the regularity assumption $v=L_\knl p$ for some $p\in L^2(\Omega)$ and that the Bernstein inequality in \eqref{eqn_BIneqC3} holds for a rate function $R_I(N)\to 0$ as $N\to \infty$.  
For a convergent Galerkin approximation scheme,  there exists a  constants $C_I,C_2>0$ such that we  have the single step error bound
\begin{align}
     \Vert v&-v_N^\delta \Vert _{L^\infty(\Omega)}\leq  \Vert v-v_N^\delta \Vert _{\sH} \nonumber \\
     &\leq C_1 \frac{\delta}{R_I(N)} \nonumber \\
     & \quad + C_2 \sup_{\xi\in \Omega} \sqrt{\knl(\xi,\xi)-\knl_{\Xi_N}^\mathrm{T}(x)\KK_N^{-1}\knl_{\Xi_N}(x)} \Vert L_\knl^{-1}v \Vert _{L^2(\Omega)} \nonumber \\
    &=C_I\frac{\delta}{R_I(N)} + C_2\sup_{\xi\in \Omega} \mathcal{P}_{H_N}(\xi)  \Vert L^{-1}_\knl v \Vert _{L^2(\Omega)},
\end{align}
with $\delta= \Vert b-b^\delta\Vert_{L^2(\Omega)}$, $\knl_{\Xi_N}(x)=\{\knl_{\xi_1}(x),\ldots,\knl_{\xi_N}(x)\}^\mathrm{T}\in \RR^N$,  $\KK_N\triangleq [\knl(\xi_i,\xi_j)]\in \RR^{N\times N}$, and $\mathcal{P}_{H_N}$ the power function of the the subspace $H_N$ in the RKHS  $\sH$. 
\end{theorem}
\begin{proof}
    The proof follows from the comments preceding the theorem. 
\end{proof}

\paragraph{\underline{Single Step Error Bounds without Bernstein Inequalities}}
The analysis above is based on the existence of Bernstein inequalities for the pair of spaces $\calW\hookrightarrow L^2(\Omega)$. Examples or cases when such inequalities hold for RKHSs is a active area of research, see the discussions in \cite{griebel2015multiscale,narcowich2012}. However, even if the kernels $\knl$ used to construct $\sH$ are ``standard ones,'' like the Sobolev-Matern kernels in \cite{narcowich2012}, the tailored kernels $\wnl$ that define $\calW$ need not be.  As emphasized in \cite{griebel2015multiscale}, an understanding of general settings when such inverse inequalities hold continues to pose a challenge to research in statistical and machine learning theory via RKHS. 

As an alternative to the explicit use of Bernstein inequalities, there are some  standard approaches, again popularized in the theory of inverse problems.  For example, these are  described in Theorem 4.2 of \cite{kress1989linear}, which show that the bound on $\|G_N\|$ in Theorem \ref{theorem_method_C2} can be replaced by 
\begin{align}
    \|G_N\| \leq C(N) \triangleq C \frac{1}{\sigma_{\min}(\AA_N\KK_N^{-1})\sigma_{\min}(\AA_N)} \label{eqn_no_bernstein}
\end{align}
via  a suitable redefinition of the interpretation of  the uncertainty to a coordinate-wise description. Note that, even without knowledge of the admittedly abstract Bernstein inequalities, the expression above can always be calculated for a fixed choice of the centers $\Xi_N$ by calculating singular valued decompositions of the matrices in the coordinate realizations.

\subsection{Approximation of the Actor Update Step}
\label{section_actorapproximation}
In the last section, two basic functional analytic settings are studied for the approximation of the critic update step in the PI iteration.
In this section, we make a systematic study of the actor update step using the RKHSs $\calU\triangleq \calU_1 \times \cdots \times \calU_m$ introduced in Section \ref{section_feature_Uspace}. The actor update step seeks to solve the operator equation 
$
\mu=Bv\in \calU$, 
where we just write $\mu\triangleq\mu^{\ell}\in \calU$ and $v\triangleq v^{\ell}\in \sH$ in this section for the actor update in \eqref{eqn_ideal2}. Because we define $B:\sH\to \calU$, we have the component-wise actor update
$
\mu_i=B_iv \in \calU_i
$ 
where $B_i\in \calL(V,\calU_i)$ for $1\leq i \leq m$. 

Via the construction in Section \ref{section_feature_Uspace}, the subspace $\calU_i=\Ran{B_i}$ is the RKHS that is generated by the kernel $\unl_i$ given in \eqref{eqn_defn_unl_i}. We introduce the set of $M$ centers $
\Omega_M =\{\omega_{1},\ldots, \omega_{M}\}\subset \Omega$, 
the associated finite-dimensional space of  approximations
\begin{align}
    U_{i,M} &\triangleq \text{span}\{\unl_{i,\omega_k}(\cdot)\triangleq \unl_i(\cdot,\omega_{k}) :  \omega_k\in \Omega_M, 1\leq k \leq M\} \nonumber \\
    &\subset \Omega,    
\end{align}
 the Grammian matrix $\UU_{i,M}\triangleq [\unl_i(\omega_j,\omega_k)]\in \RR^{M\times M}$,  the $\calU_i$-orthogonal projection $\Pi_{U_{i,M}}:\calU_i\to U_{i,M}$ onto $U_{i,M}$, and $\mathcal{P}_{U_{i,M}}$ the power function of the subspace $U_{i,M}$ in $\calU_i$.   For any $\mu_i\in \calU_i$, the coordinates $\alpha_M\triangleq \{\alpha_1,\ldots,\alpha_M\}^{\rm T}\in \RR^{M\times 1}$ of the realization 
 $
 \Pi_{U_{i,M}}\mu_i\triangleq \sum_{k=1,\ldots,M} \alpha_k \unl_{i,k}
 $
 satisfy the matrix equation 
 \begin{align}
     \sum_{k=1}^M\unl_{i}(\omega_j,\omega_k) \alpha_k =& \mu_i(\omega_j) \quad \text{ for all } 1\leq j\leq M
 \end{align}
 or, equivalently,
 \begin{align}
     \UU_{i,M} \alpha_M&=\MM_M,
 \end{align}
 where $\MM_M\triangleq \left[\mu_{i}(\omega_1),\ldots, \mu_i(\omega_M)\right]^{\rm T}\in \RR^M$.
With this definition of the approximation subspaces for the actor step, we obtain the following result. 

\begin{theorem}
    \label{theorem_actorapprox}
Let $\mu\triangleq \mu^\ell=\left[\mu_1,\dots,\mu_m\right]^{\rm{T}}\in \calU$ be the solution of the abstract operator equation in the actor step in \eqref{eqn_ideal2}. Define the approximation of the actor step $\mu_{i,(M,N)}\triangleq \mu^\ell_{i,(M,N)}$ to be given by  
\begin{align}
    \mu_{i,(M,N)}\triangleq \Pi_{U_{i,M}}Bv^\delta_N, 
\end{align}
where $v^\delta_N$ is the Galerkin approximation of Method $(C1)$ or Method $(C2)$ described in  Theorems \ref{theorem_method C1} or \ref{theorem_method_C2}, respectively. 
For convergent Galerkin approximations in the critic and actor steps of the PI algorithm, 
there exist three constants $C(N),C_u,C_v>0$ such that 
\begin{align}
     \Vert \mu_i&-\mu_{i,(M,N)} \Vert _{L^\infty(\Omega)} \nonumber\\
     &\leq C_u \sup_{\xi\in \Omega}\mathcal{P}_{U_{i,M}}(\xi)\Vert \mu_i \Vert _{\calU_i} + C_v  \Vert v-v^\delta_N \Vert _{\sH} \nonumber \\
    &\leq C(N)\delta + C_v\sup_{\xi\in \Omega} \mathcal{P}_{H_N}(\xi)  \Vert L^{-1}_\knl v \Vert _{L^2(\Omega)} \nonumber \\
    & \quad + C_u \sup_{\xi\in \Omega}\mathcal{P}_{U_{i,M}}(\xi)) \Vert \mu_i \Vert _{\calU_i}.
\end{align}
\end{theorem}
\begin{proof}
    We apply the triangle inequality to obtain
    \begin{align}
        |E_x( \mu_i&-\mu_{i,(M,N)})| \nonumber \\ &\leq  \left|E_x(I-\Pi_{U_{i,M}})\mu_i \right| + \left|E_x\Pi_{U_{i,M}}B_i(v-v_N^\delta ) \right|, \nonumber \\
        &\leq \sup_{\xi\in \Omega} \mathcal{P}_{U_{i,M}}(\xi)  \Vert \mu_i \Vert _{\calU_i} + \bar{\unl}_i \Vert B_i \Vert   \Vert v-v_N^\delta \Vert _{\sH}.
    \end{align}
    In this inequality, we have used the fact that, owing to  the regularity assumption $\knl\in C^{2s}(\Omega\times \Omega)$,  as well as the fact that  $\Omega$ is compact, each reproducing kernel $\unl_i$ is bounded uniformly on the diagonal. There  is a constant $\bar{\unl}_i>0$ such that  
    $\unl_i(x,x)\leq \bar{\unl}_i^2$ for all  $x\in \Omega$. 
    This implies that the family of evaluation functions $\{E_x\}_{x\in \Omega}$, viewed as mappings $E_x:\calU_i\to \RR$,  have uniformly bounded operator norms, with 
    $\|E_x\|_{\calL(\calU_i,\RR)}\leq \bar{\unl}_i$ for all $x\in \Omega, 1\leq i\leq M$. 
    For Method  $(C1)$, the coefficient $C(N)$ can be taken as $C_1$ in Theorem \ref{theorem_method C1}, while for Method $(C2)$,
    we can take $C(N)=C_I/R_I(N)$ as in Theorem \ref{theorem_method_C2}.
    For the $(C2)$ Method, in the absence of Bernstein inequalities, $C(N)$ can be taken as in \eqref{eqn_no_bernstein}. 
\end{proof}

\section{Numerical Examples}
In this section, we describe empirical numerical results obtained when the tailored RKHS  described in this paper are used to achieve optimal control approximations via policy iteration.
We study a nonlinear system of the form in \eqref{eqn_sys}, which has been used as a testbed by some of the authors in \cite{vamvoudakis2010online}.
We let $n = 2$, $m = 1$,
$
f(x) =[ -x_1 + x_2 ,
-0.5 x_1 - 0.5 x_2  ( 1 -  ( \cos(2x_1) + 2  )^2  )]^{\rm T}
$
and
$g(x) =  [0,\cos(2x_1) + 2]^{\rm T}.$ 
We define the standard cost function of the linear-quadratic regulator problem in \eqref{eqn_J_lqr} with $ R = 1 $ and $ Q = I_2 $, where $ I_2 $ denotes the $2 \times 2$ identity matrix.
For this cost function, the optimal value function is known in closed form and is given by
$
V^\star(x) = 0.5 x_1^2 + x_2^2,
$
and the optimal control policy is
$
\mu^\star(x) = -(\cos(2x_1) + 2)x_2.
$

Studies of the performance of the PI algorithm were carried out for the Gaussian or exponential kernel, as well as  for  three Sobolev-Mat\'{e}rn kernel functions having various smoothness indices $k$. The exponential kernel is written as $\knl(x_1,x_2)=\eta({\|x_1-x_2\|_2})$ in terms of the radial function $\eta$ given by $
\eta (r) = e^{-\alpha r^2}$, 
were $\alpha$ is known as the hyperparameter, which is selected to be $1$ in this section.
The Sobolev-Mat\'{e}rn kernel having integer smoothness index $k>0$ is defined in \cite[Ch. 3]{RKHS_MRAC_book}.
The Sobolev-Mat\'{e}rn kernel of smoothness $k$  and defined over  subsets $\Omega\subset \RR^n$  is written in terms of a radial function as $\knl(x_1,x_2)\triangleq \eta_k(\|x_1-x_2\|_2)$, where 
\begin{align}
    \eta_k (r) = \frac{2\pi^{n/2}}{\Gamma(k)}{K}_{k-n/2}(r/2)^{k-n/2}.    
\end{align}

The primary purpose of these numerical simulations is to study the qualitative relationship between the single-step error bounds derived in the paper and the overall performance of the PI algorithm after completing all the steps.
We emphasize that,  so far, we have not studied the rates of convergence of the overall PI algorithm as the number of iteration steps $\ell\geq 0$ increases.
For this reason, simple choices are made for the location of the centers.
Specifically, the centers are uniformly distributed in the domain $\Omega = [-1,1] \times [-1,1]$. 

In the simulation results summarized below, we compared the overall PI convergence rates, the performance of the learned approximate optimal controllers, and the computational costs of the PI iteration using the Methods $(C1)$ and $(C2)$ for the critic
As is usual in many studies of approximations in RKHS \cite{wendland,schaback94},
we chose to present the comparisons of convergence in terms of the fill distance 
\begin{align}
    h_{\Xi_N,\Omega}\triangleq \sup_{\xi\in \Omega}\min_{\xi_i\in \Xi_N}\|\xi-\xi_i\|_{\RR^n}.    
\end{align}
of the centers $\Xi_N\subset \Omega$ that are used to define the approximants $H_N$ of the value function.
Also, while the approximation of the actor step described in Section \ref{section_actorapproximation} allows for a different set of centers for approximating  $\mu$, for simplicity,  we use the same set, choosing $M=N$ and  $\Omega_M\triangleq \Xi_N$. 

%
%
\subsection{The Comparison of  Convergence Across Iterations}
  In general, Algorithm \ref{algorigm_PI_2} can be implemented when the dimension $N(\ell)$ of approximations varies over iterations.
  It is likewise simple to replace the termination criteria in terms of the maximum number of iterations $\bar{\ell}$ by one defined by a desired convergence tolerance $\epsilon>0$.
  We choose a very simple setting for the benchmark simulations summarized in this initial paper, just fixing the final number of iterations $\bar{\ell}$.  

\begin{figure*}[ht]
    \centering
    \includegraphics[width=1\linewidth]{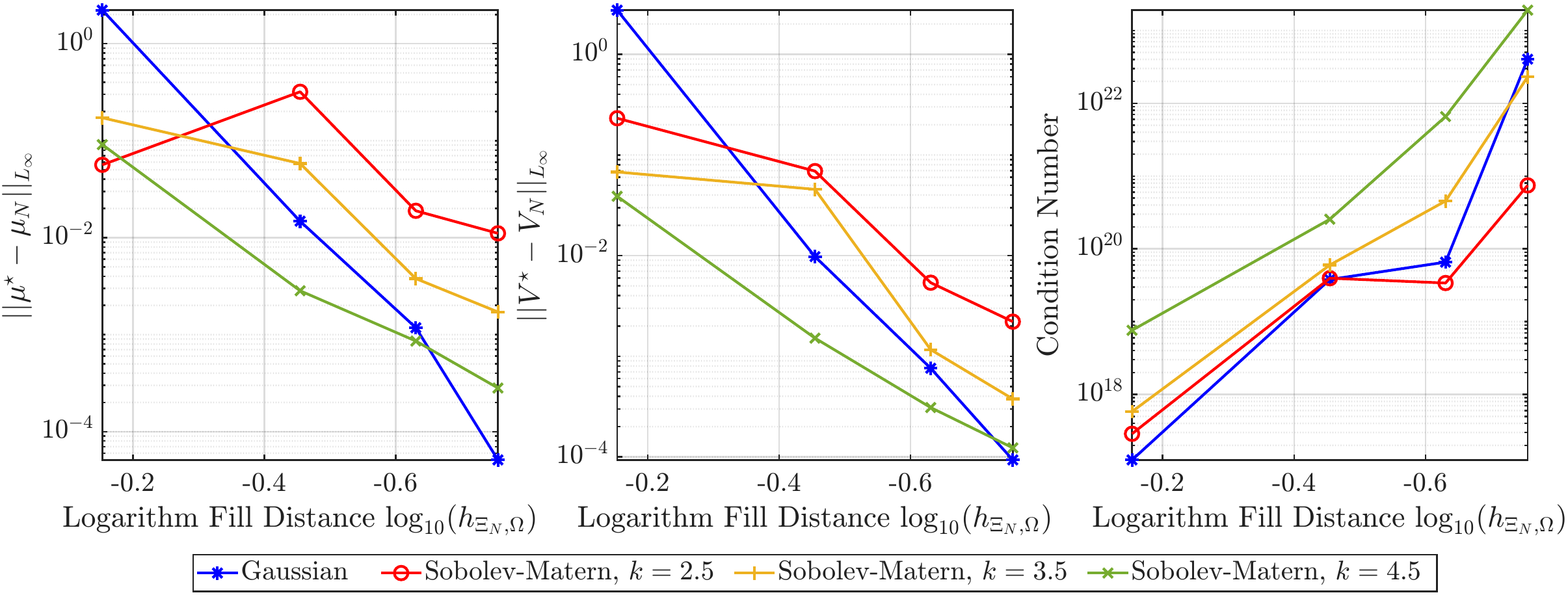}
    \caption{Convergence results of PI  using the (C1) Method. In these figures the maximum number of iterations  $\bar{\ell}$ in Algorithm \ref{algorigm_PI_2} is set to $\bar{\ell}=40$}
    \label{figure_C1_result}
\end{figure*}

  In Figure \ref{figure_C1_result}, the performance of the approximately optimal controller and associated value function obtained in the final step of the PI algorithm, for different choices of the reproducing kernel, is plotted versus the fill distance of a uniform set of centers.
 The number of centers used across the different iterations $\ell\leq \bar{\ell}=40$ of the PI algorithm is held fixed in these benchmark tests.  The maximum number of iterations $\bar{\ell}$, at least intuitively, is chosen to be quite large, with the goal of conveying what might be thought of as the practical limit of achievable performance of the PI method. 
 
As expected, using the (C1) Method, the overall or final errors achieved decrease as $N \to \infty$, for the fixed large number of maximum number of $\bar{\ell}=40$ iterations in this comparison.
In the figures, the increase in $N$   corresponds to the fill distance $h_{\Xi_N,\Omega}\to 0$. 
It is well-known that smoother kernels lead to faster rates of convergence in the usual setting when learning theory is used to construct interpolants for regression problems (see Table 11.1 of \cite{wendland} or \cite{schaback94}), and similar phenomenology is apparent for Method  $(C1)$ for our problem of formulating the  PI algorithm in RKHS spaces.   
The Gaussian kernel is the smoothest kernel included in the comparison, and for problems of learning theory or interpolation for regression, it is known to have a spectral order of convergence, see the discussion in \cite[Sec. 11.4]{wendland}.
However, it is sometimes described as generating a very small native space, one that is contained in $C^\infty(\Omega)$.
In our applications, this means that the collection of ODE systems for which $V^*\in \sH$,  and for which optimality can be guaranteed,  is correspondingly small, relatively speaking.
It exhibits the fastest rate of PI convergence over iterates in Figure 2, which is consistent with the expectations from applications in regression and the bounds in Table 11.1 of \cite{wendland}.
Also, the error associated with the controller or value function is smaller when the Sobolev-Mat\'{e}rn kernel function is smoother. However, it is not apparent that the \textit{rate of convergence} of the error for the Sobolev-Mat\'{e}rn kernel is faster as the smoothness index $k$ increases.
This is a notable qualitative difference compared to standard results for regression via learning theory in RKHS.  This observation  could be an artifact of the rather simple assumptions of fixed dimensionality of approximations for each individual step, or it could result from  the fixed maximum number of iterations $\bar{\ell}=40$.
Further refined studies are warranted on how changing the dimension and number of iterations affects performance.
This would arise naturally in the study of how to adaptively choose the individual basis functions, or how to adaptively choose the number of bases, a topic discussed further below. 

The plot on the right represents the condition number of ``a typical''  matrix $\AA_N$  during the computation. The condition number is calculated based on the matrix $\AA_N\triangleq \AA_{N(\bar{\ell})}$, which is the matrix in the final iteration. 
This plot shows us that the condition number grows when the number of centers increases, a well-known phenomenon studied for approximations in RKHS;
see \cite[Ch. 11]{wendland} for a background discussion. 
For Method $(C1)$, the matrix realization $\AA_N$ is in general not symmetric, and the condition number takes on exceptionally high values as $N$ increases, across all the choices of kernels tested for the  PI comparison.
Because of these high condition numbers, the truncated singular value decomposition is employed during practical computations.

\begin{figure*}[ht]
    \centering
\includegraphics[width=1\linewidth]{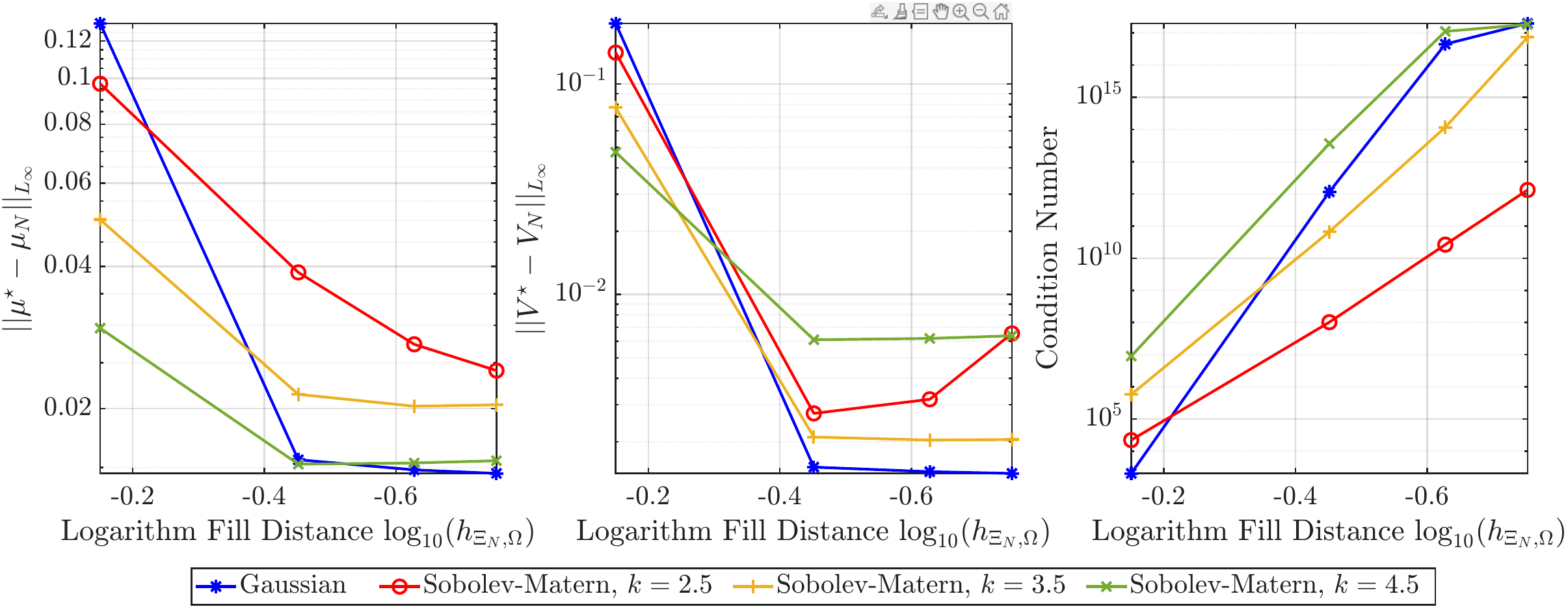}
    \caption{Convergence results of PI  using Method($C2$). In these figures, the target number of iterations  $\bar{\ell}$  in Algorithm \ref{algorigm_PI_2} is set to $\bar{\ell}=40$}
    \label{figure_C2_results}
\end{figure*}

The study of the PI algorithm using Method ($C2$) further supports some observations made in \cite{bouland2023}. However, it is important to note that this paper extends those conclusions and makes observations for the complete PI algorithm, including both critic and actor updates working together.

We used an \textit{ad hoc} approach to selecting quadratures.
All the testbed studies were conducted with an increasing number of quadrature points, and the results are reported for several quadrature points high enough that no discernible differences are apparent in the final approximately optimal controller. Again, a careful characterization of the effect of quadrature error remains an important topic for future research. However, under these conditions, the same phenomenology is apparent in this case as in the case above for Method $(C1)$: smoother kernels yield smaller errors in the final iterate of the PI algorithm.

Note that the \textit{rate of convergence} of the error in the final estimate of the value function when using Method $(C2)$ increases with smoother kernels. This stands in contrast to the situation noted above using Method $(C1)$. 

As expected from the stepwise error bound for Method $(C2)$, the contribution from the first term in the error bound,
$\dfrac{1}{\sigma_\textrm{min}(\AA_N\mathbb{\KK}^{-1}{N}) \sigma_{\min}(\AA_N)}\delta$,
becomes dominant as $N \to \infty$.
The plot on the right illustrates the corresponding condition number for $\AA_N$ when using Method $(C2)$.
Again, while smoother kernels lead to smaller errors in the optimal controller, their associated condition numbers tend to grow faster.
However, for Method $(C2)$, the matrix realizations $\AA_N$ and $\KK_N$ are symmetric, unlike the situation with Method $(C1)$.
Perhaps most importantly, typical condition numbers can be significantly smaller for kernels with lower smoothness, which is a crucial difference compared to the results obtained using Method $(C1)$. 

\subsection{Performance Comparison}

\begin{figure*}[t]
    \centering
    \includegraphics[width=1\linewidth]{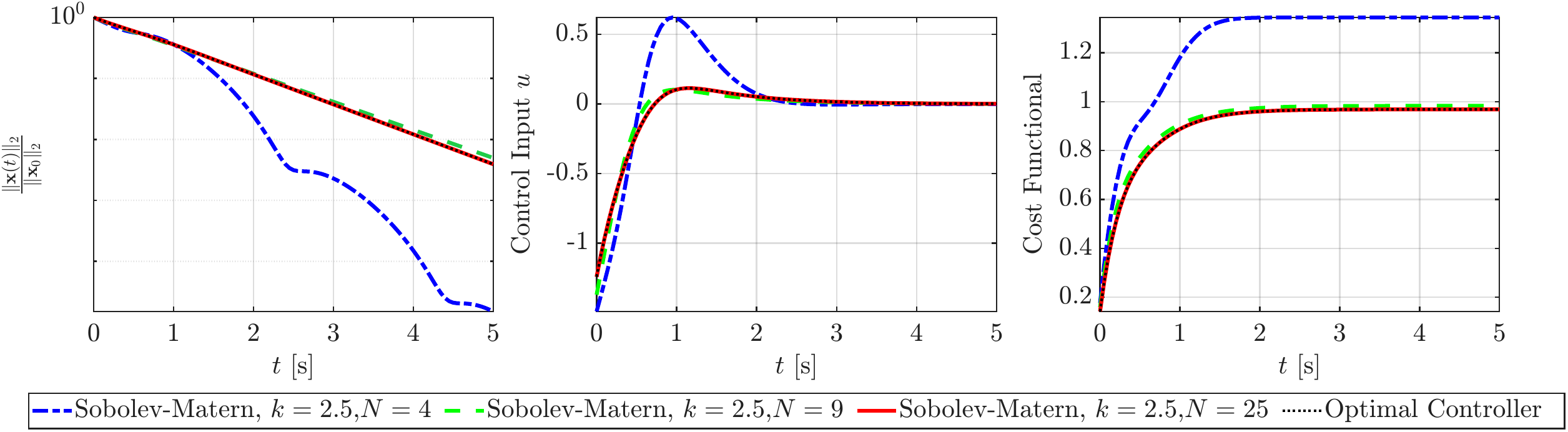}
    \caption{Transient simulation with learned control from Method $(C1)$}
    \label{C1_test}
\end{figure*}
\begin{figure*}[t]
    \centering
    \includegraphics[keepaspectratio=true,width=1\linewidth]{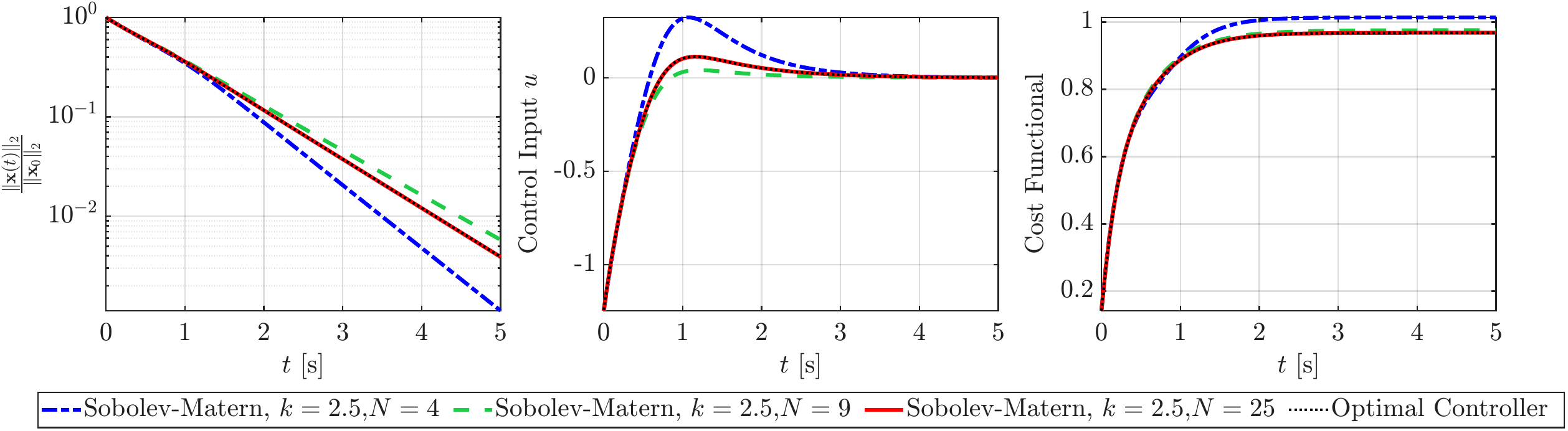}
    \caption{Transient simulation with learned control from Method $(C2)$}
    \label{figure_C2_TimeRes}
\end{figure*}

The previous section studied the asymptotics of the errors in the controllers and value functions generated by the PI algorithm as the dimension increases, when the total number of iterations $\bar{\ell}$ is held constant.
However, in practice, it is important to assess what constitutes a practical choice of $N$.
In these initial comparisons, we only observe that often choosing $N$ to be small works sufficiently well. 
In Figures \ref{C1_test} and \ref{figure_C2_TimeRes}, the transient simulations of three different controllers for both Method $(C1)$ and Method $(C2)$ are presented for demonstration purposes, with initial condition $x_0 = [0.9,0.7]$. 
The three controllers' response, control effort, and cost function are compared to the optimal controller over time.
As expected, the smoothest kernel with a higher number of centers gives the best approximation of the optimal controller for both methods.  
However, with a smaller number of centers, both methods can still generate reliable controllers that can also drive the system states to the origin.
Still,  the controllers are less energy efficient.
While the error bounds derived in this paper are a good start, methods for adaptively choosing a practical, but adequate, number of bases $N$ are still very much needed.
We comment further on this fact in the conclusions and recommendations for future work. 


\begin{table*}[h]
\centering
\footnotesize 
\setlength{\tabcolsep}{5pt} 
\renewcommand{\arraystretch}{1.2} 
\begin{tabular}{|c|cccccccc|}
\hline
\multirow{3}{*}{Kernel Functions} & \multicolumn{8}{c|}{Number of Centers} \\ \cline{2-9} 
                                  & \multicolumn{2}{c|}{9}  & \multicolumn{2}{c|}{25} & \multicolumn{2}{c|}{49} & \multicolumn{2}{c|}{81} \\ \cline{2-9} 
                                  & \multicolumn{1}{c|}{$C1$} & \multicolumn{1}{c|}{$C2$} & \multicolumn{1}{c|}{$C1$} & \multicolumn{1}{c|}{$C2$} & \multicolumn{1}{c|}{$C1$} & \multicolumn{1}{c|}{$C2$} & \multicolumn{1}{c|}{$C1$} & $C2$ \\ \hline
Gaussian                          & \multicolumn{1}{c|}{0.011} & \multicolumn{1}{c|}{0.394} & \multicolumn{1}{c|}{0.017} & \multicolumn{1}{c|}{0.429} & \multicolumn{1}{c|}{0.062} & \multicolumn{1}{c|}{0.435} & \multicolumn{1}{c|}{0.124} & 0.438 \\ \hline
Sobolev-Matern, $k = 2.5$         & \multicolumn{1}{c|}{0.007} & \multicolumn{1}{c|}{0.385} & \multicolumn{1}{c|}{0.036} & \multicolumn{1}{c|}{0.648} & \multicolumn{1}{c|}{0.072} & \multicolumn{1}{c|}{0.679} & \multicolumn{1}{c|}{0.151} & 0.714 \\ \hline
Sobolev-Matern, $k = 3.5$         & \multicolumn{1}{c|}{0.009} & \multicolumn{1}{c|}{0.549} & \multicolumn{1}{c|}{0.025} & \multicolumn{1}{c|}{1.172} & \multicolumn{1}{c|}{0.074} & \multicolumn{1}{c|}{1.208} & \multicolumn{1}{c|}{0.156} & 1.194 \\ \hline
Sobolev-Matern, $k = 4.5$         & \multicolumn{1}{c|}{0.009} & \multicolumn{1}{c|}{0.810} & \multicolumn{1}{c|}{0.023} & \multicolumn{1}{c|}{1.816} & \multicolumn{1}{c|}{0.074} & \multicolumn{1}{c|}{1.906} & \multicolumn{1}{c|}{0.162} & 1.817 \\ \hline
\end{tabular}
\caption{Computational time in seconds for 40 iterations of Methods $(C1)$ and $(C2)$ methods in a PC equipped with an Intel Core i7-11700 processor}
\label{table_cputime}
\end{table*}
 
The studies to date have compared the performance of the overall PI algorithm with various types of kernels.
We have demonstrated that assessing the advantages and disadvantages of these alternatives in relative performance can be subtle.
It is also important to observe that the differences in computational work between the Methods $(C1)$ and $(C2)$ can be substantial.
As noted in the definition of Method $(C2)$, practical computations require the use of quadrature approximations.
The cost of using quadratures in Method $(C2)$ dominates the computational work of this method, and it is significantly more costly overall to implement than Method $(C1)$. 
In Table \ref{table_cputime}, we illustrate the computational time for a PC equipped with an Intel Core i7-11700 processor, comparing the performance of both methods.
The difference in computational cost varies between $O(10)$-$O(100)$ over the range of the dimension $N$.

\section{Conclusion}
This paper introduced a novel framework for tailoring the kernel and for approximating of nonlinear optimal control problems using policy iteration.
By constructing RKHSs that are intrinsically defined with respect to the system dynamics, we derived rigorous stepwise performance bounds that hold for each iteration of the PI algorithm.
Our analysis extended our recent work \cite{bouland2023} to include the actor update, completing a cohesive kernel-based approach to reinforcement learning in control. We also introduced two Galerkin approximation methods for the critic update, each with distinct trade-offs in computational cost and numerical sensitivity.
In short, Method $(C1)$ is less computationally costly, but Method $(C2)$ exhibits better numerical conditioning properties in these initial implementations and tests.  

While, to the authors' knowledge, these are the first such closed-form, stepwise error bounds for the PI algorithm including numerical perturbation terms, they only constitute a first step in the rigorous study of how methods of statistical and machine learning theory can improve our understanding of some classical algorithms in nonlinear optimal control.
These results bridge optimal control theory and statistical learning, potentially paving the way for further study and development of adaptive, structure-aware control algorithms with provable performance guarantees.
However, several important questions remain open.
Examples of these questions include the development of associated methods for the adaptive selection of bases, the development of practical termination criteria during basis enrichment or pruning, and the further study of the contribution of quadrature errors to PI performance.  

\section{Appendices}

\subsection{Elements of Reproducing Kernel Hilbert Spaces}\label{sec:rkhs}
For a systematic presentation of RKHS theory using the same notation as in this paper, see \cite[Ch. 3]{RKHS_MRAC_book}. 
Let the kernel $\knl:\Omega \times \Omega \to \RR$ define the RKHS $\sH$ of real-valued functions over $\Omega \subseteq \RR^n$.
As in \cite{guella2022operator,zhou2008derivative},
in this paper, we assume that $\Omega$ is compact,
simply connected, and 
$\overline{\Omega}=\overline{\left(\mathring{\Omega}\right)}$.
This assumption excludes sets containing isolated points and
allows extending derivatives from $\mathring{\Omega}$ to all of $\Omega$.
Following the standard notation in the study of ODEs and PDEs,
we denote a mixed partial derivative of a real-valued function
$f : \Omega \subseteq \RR^n\to \RR$ for a non-negative integer multiindex $\alpha\triangleq (\alpha_1,\ldots,\alpha_n)\in (\ZZ^+_0)^n$ by 
    $D^\alpha f(x)\triangleq {\partial^{|\alpha|} f(x)}/{\partial x^{\alpha_1}_1 \cdots \partial x^{\alpha_n}_n}$,  
where $|\alpha| \triangleq \sum_{i=1}^n \alpha_i$ denotes the {length of} $\alpha$.
For any non-negative integer $m$,
we let $C^m(\mathring{\Omega})$ denote the class of $m$-times continuously differentiable functions with norm
     $\Vert f \Vert _{C^m(\mathring{\Omega})} \triangleq \sum_{0\leq |\alpha|\leq m}  \Vert D^\alpha f \Vert _{C(\mathring{\Omega})}$
for any  $f \in C^m(\mathring{\Omega})$. 
The space $C^{m}(\Omega)$ over the compact set $\Omega$ contains all the  extensions to the boundary of all the  functions in $C^{m}(\mathring{\Omega})$,
equipped with the norm 
     $\Vert f \Vert _{C^m(\Omega)}\triangleq \left  \Vert f|_{\mathring{\Omega}}\right  \Vert _{C^m(\mathring{\Omega})}$ for any  $f \in C^m(\Omega)$. 
We slightly generalize our notation when discussing the derivatives of kernel functions $\knl : \Omega\times \Omega \to \RR$.
Specifically,
we denote by $(\alpha,\beta)$ a pair of non-negative integer multi-indices, each of length $n$, and
write $D^{(\alpha,\beta)}\knl(x,y)$ as 
\begin{align}
   D^{(\alpha,\beta)}\knl(x,y) &\triangleq 
   \frac{\partial^{|\alpha|+|\beta|} \knl(x,y)}{\partial x^{\alpha_1}_1 \cdots \partial x^{\alpha_n}_n \partial y^{\beta_1}_1 \cdots \partial y^{\beta_n}_n}.   
\end{align}

\begin{theorem}[{\cite[Th. 1]{zhou2008derivative}}]\label{theorem_zhou}
    Let the kernel $\knl : \Omega\times \Omega \to \RR$ generate the RKHS $\sH$ and
    assume that 
    $\knl\in C^{2s}(\Omega\times \Omega)$ with $s \geq 1$.
    Then,
    \begin{enumerate}
        \item For any $x\in \Omega$ and $|\alpha|\leq s$, it holds that
        \begin{align}
            D^{(\alpha,0)}\knl(x,\cdot)\in \sH;
        \end{align}
        \item For all $x\in \Omega$ and $f\in \sH$,
        it holds  
        \begin{align}
            D^\alpha f(x)=\left (D^{(\alpha,0)}\knl(x,\cdot),f \right)_\sH;
        \end{align}
        \item The inclusion map $\mathcal{I}:\sH \hookrightarrow C^s(\Omega)$ is continuous, compact, and well-defined, That is, 
        for any compact subset of functions $\mathcal{C}\subset \sH$,  the collection of functions $\mathcal{I}(\mathcal{C})\subset C^s(\Omega)$ is compact. 
    \end{enumerate}
\end{theorem}

We close this quick review of RKHSs
by noting that in this paper,
we always construct approximations by selecting some subset of $N$ centers $\Xi_N\triangleq \{\xi_1,\cdots,\xi_N\}\subset \Omega$ and defining the associated finite dimensional subspace 
\begin{align}
    H_N\triangleq \text{span}\left \{ \knl_{\xi_i} \in \Xi_N :  1\leq i\leq N\right \}\subseteq \sH. \label{eqn_finite_dimensional_space}
\end{align}
Given our standing assumption that the kernel $\knl$ is strictly positive, the dimension of $H_N$ is $N$.
We denote by $\Pi_N:\sH\to H_N$  the $\sH$-orthogonal projection of $\sH$ onto $H_N$. It is well-known that the orthogonal projection $\Pi_Nf$ interpolates the function $f\in \sH$ at the centers, and we always have the pointwise error 
\begin{align}
    |E_x(I-\Pi_N)f|&\leq \Pwr_N(x)\|(I-\Pi_N)f\|_\sH \nonumber \\
    &\leq \Pwr_N(x)\|f\|_\sH \,\,\, \text{ for all } x\in \Omega, f\in \sH, 
\end{align}
where
\begin{align}
    \Pwr_N(x)\triangleq \sqrt{\knl(x,x)-\knl_N(x,x)} \quad \quad \text{ for all } x\in \Omega,
\end{align}
denotes the power function, and
$\knl_N:\Omega \times \Omega \to \RR$ is the reproducing kernel of the subspace $H_N$ defined by 
\begin{align}
    \knl_N(x_1,x_2)\triangleq \knl_{\Xi_N}^{\rm T}(x_1)\KK_N^{-1} \knl_{\Xi_N}(x_2),
\end{align}
with $\KK_N$ denoting the Grammian matrix of the centers and the column vector of basis functions  $\knl_{\Xi_N}(x) \triangleq \left[\knl_{\xi_1}(x),\ldots, \knl_{\xi_N}(x)\right]^{\rm T}\in \RR^N$ for all $x\in \Omega$. 

\subsection{Theory  of Galerkin Approximations}\label{section_galerkinapproximations}
 This paper uses several properties of Galerkin approximations of operator equations, and this section summarizes the needed theory.
 An in-depth and comprehensive discussion of the topics below can be found in the references on inverse problems and their approximations. See \cite{engl1996regularization,kress1989linear,kirsch2011introduction}, which provide full discussions and contain the original versions of the theorems we summarize in this appendix. 

 We consider the noise-free operator equation 
 $Tx=y$, 
where $T\in \calL(X,Y)$ is an injective, bounded, linear,  operator acting between the normed vector spaces $X$ and $Y$,  and its associated noisy operator equation is defined to be $Tx^\delta = y^\delta$, for some $y^\delta$ such that $\|y-y^\delta\|_Y\leq \delta$.    
In a general Petrov-Galerkin approximation, we  introduce the finite-dimensional subspaces 
    $X_N\triangleq \text{span}\{\phi_1,\ldots,\phi_N\}\subseteq X$,  $ 
    Y_N\triangleq \text{span}
    \{\psi_1,\ldots,\psi_N\}\subseteq Y$, and denote by $\{P_N\}_{N\geq N_0}$ a family of bounded linear projection  operators $P_N:Y\to Y_N$,  with $P_N\in \calL(Y)$. When we say that $P_N$ is a projection operator, this means it is idempotent, with  $P_N^2=P_N$ for each $N$.  We define the Petrov-Galerkin approximation  $x_N\in X_N$ generated by $y\in Y$ to be given by 
$x_N \triangleq  {\left (P_N T|_{X_N} \right)^{-1}P_N}\triangleq { G_N}y$, 
where $G_N$ is the Galerkin operator. The collection of subspaces and projection operators $\{(X_N,Y_N,P_N)\}_{N\geq N_0}$ defines the Petrov-Galerkin approximation  scheme. 

Since $P_N:Y\to Y_N\subset Y$ is a bounded linear operator $P_N\in \calL(Y)$, it has the general representation 
$
P_Ny=\sum_{i=1,\ldots,N} \lambda_{N,i}(y)\psi_i$,
where $\lambda_{N,i}:Y\to Y_N$ are uniquely defined  bounded linear coordinate functionals that define $P_N$. The matrix representations of the Galerkin approximations of the noise-free operator equation then take the form 
$
\TT_N \Theta_N = Y_N$ 
where $\TT_{N}\triangleq [\lambda_{N,i}(A\phi_j)]\in \RR^{N\times N}$, $Y_N\triangleq \{\lambda_{N,1}(y),\ldots \lambda_{N,N}(y)\}^{\rm T} \in \RR^N$, $\Theta_N\triangleq \{\theta_1,\ldots, \theta_N\}^{\rm T}\in \RR^N$, and $x_N\triangleq \sum_{j=1}^M \theta_j\phi_j$. 

If $Y$ is a Hilbert space, then we can choose $P_N\triangleq \Pi_{Y_N}$ to be the orthogonal projection of $Y$ onto $Y_N$.
In this special case, the equations that define the Galerkin approximation $x_N\in X_N$ of the solution $x$ in the noise-free equation can be written in terms of inner products
$\left ( Tx_N-y,\psi_i\right)_Y = 0$ for all $1\leq i\leq N$. 
The coordinate realizations of these equations  consequently seek  to find the coefficients $\{\theta_j\}_{j=1}^N\in \RR^N$ in $x_N(\cdot)=\sum_{j=1}^N \phi_i(\cdot)\theta_i$ such that 
$\sum_{j=1,\ldots,N} \left (\psi_i,T\phi_j \right)_Y\theta_j  = \left (\psi_i,y \right)_Y$ for  $i=1, \dots, N$. This is written succinctly as 
$\TT_N \Theta_N=Y_N$ 
with $\TT_N\triangleq [(\psi_i,T\phi_j)_Y]\in \RR^{N\times N}$ and
$Y_N\triangleq \left[(\psi_1,y)_Y,\ldots, (\psi_N,y)_Y\right]^{\rm T}\in \RR^N$. A similar representation holds for the Galerkin approximation $x^\delta_N\in X_N$ of the solution $x^\delta$ of the noisy operator equation. 

There are some standard properties of the Petrov-Galerkin operator $G_N$ that play an important part of the analysis in this paper. We say that the family of subspaces $\{X_N\}_{N\geq N_0},\{Y_N\}_{N\geq 0}$ and operators $\{P_N\}_{N\geq N_0}$ define a convergent Galerkin scheme if $G_N$ is well-defined for each $N\geq N_0$,  and for each $y\in \Ran{A}$ we have $\lim_{N\to \infty}G_N y \to x\triangleq T^{-1}y$ as  $N\to \infty.$

\begin{theorem}[\cite{kirsch2011introduction,kress1989linear}]\label{theorem_propsGalerkin}
Suppose that $T$ is an injective bounded linear operator and that $G_N\triangleq \left (P_NT|_{X_N}\right)^{-1}P_N$ is well-defined for each $N\geq N_0$. We have the following:
    \begin{enumerate}
     \item Each operator $G_NT$ is a projection on $X_N$. 
        \item The Galerkin scheme is convergent if and only if the sequence of operators $\{G_NT\}_{N\geq N_0}$ is uniformly bounded in the dimension $N$. That is, there is a constant $c>0$ such that 
         $\Vert G_N T \Vert \leq c < \infty$.
        \item If $T$ is compact and $\text{dim}(X)=\infty$, then 
         $\Vert G_N \Vert \to \infty$. 
        \item For a convergent Galerkin method, we have the following basic estimate: there exists a constant $C>0$ such that 
        \begin{align}
            \|x-x_N^\delta\|_X  \leq \|G_N\|\delta + C \inf_{g_N\in X_N}\|x-g_N\|_X, 
        \end{align}
        where $x$ is the solution of the noise-free operator equation and $x_N^\delta$ is the Galerkin approximation using noisy data as defined in the noisy operator equation. 
    \end{enumerate}
\end{theorem}
 

\bibliographystyle{IEEEtran}
\bibliography{inv,inv1,koop,tailorkernels}
 

\begin{IEEEbiographynophoto}{Shengyuan Niu} received the B.S. and M.S. degrees in Mechanical Engineering from Virginia Tech in 2020 and 2023, respectively, and is currently pursuing a Ph.D. degree in Mechanical Engineering at Virginia Tech. His research interests include optimal control, learning theory, and approximation methods using reproducing kernel Hilbert space embeddings. His current work focuses on developing approximation methods for solving the Hamilton-Jacobi-Bellman equation through reinforcement learning
\end{IEEEbiographynophoto}

\begin{IEEEbiographynophoto}{Ali Bouland} is a Ph.D. candidate in Mechanical Engineering at Virginia Tech. He received his B.S. in Mechanical Engineering from Southern Illinois University Carbondale (SIUC) and his M.S. in Mechanical Engineering from the University of Southern California (USC). His research interests include control theory and reinforcement learning, with an emphasis on developing modern methods for the analysis and control of complex dynamical systems. He is particularly interested in applications to energy systems and flow control, where his work aims to bridge advanced theoretical methods with practical engineering challenges.
\end{IEEEbiographynophoto}

\begin{IEEEbiographynophoto}{Haoran Wang} is a Ph.D. candidate of Mechanical Engineering in Virginia Tech, Blacksburg, VA, USA.
He received a B.S. degree in Mechanical Engineering from Virginia Tech in 2021 and an M.S. degree in Mathematics from the same institution in 2025.
He is currently working with Dr. Andrew J. Kurdila and Dr. Andrea L'Afflitto on combining model reference adaptive control with reproducing kernel Hilbert space embedding. 
\end{IEEEbiographynophoto}

\begin{IEEEbiographynophoto}{Filippos Fotiadis}
(Member, IEEE) was born in Thessaloniki, Greece. He received the PhD degree in Aerospace Engineering in 2024, and the MS degrees in Aerospace Engineering and Mathematics in 2022 and 2023, all from Georgia Tech. Prior to his graduate studies, he received a diploma in Electrical \& Computer Engineering from the Aristotle University of Thessaloniki. He is currently a postdoctoral researcher at the Oden Institute for Computational Engineering \& Sciences at the University of Texas at Austin.
His research interests are in the intersection of systems \& control theory, game theory, and learning, with applications to the security and resilience of cyber-physical systems.
\end{IEEEbiographynophoto}

\begin{IEEEbiographynophoto}{Andrew Kurdila}
received a B.S. degree in
engineering mechanics from the University of
Cincinnati, Cincinnati, OH, USA, in 1983, the
M.S. degree in engineering mechanics from the
University of Texas, Austin, TX, USA, in 1985,
and the Ph.D. degree in engineering science
and mechanics from Georgia Tech, Atlanta, GA,
USA, in 1989.
He is a Professor of Mechanical Engineering at Virginia
Tech, and from 2006-2022, he served as the  W. Martin Johnson Professor in the same department.
Previously, he served as a tenured Faculty Member at Texas A\&M
University, College Station, TX, USA, and the University of Florida,
Gainesville, FL, USA.
He is the author of more than 200 refereed
conference and journal articles and 3 books.
His research interests include dynamical systems theory, learning and
approximation theory, control theory, and computational mechanics.
\end{IEEEbiographynophoto}

\begin{IEEEbiography}[{\includegraphics[width=1in,height=1.25in,clip,keepaspectratio]{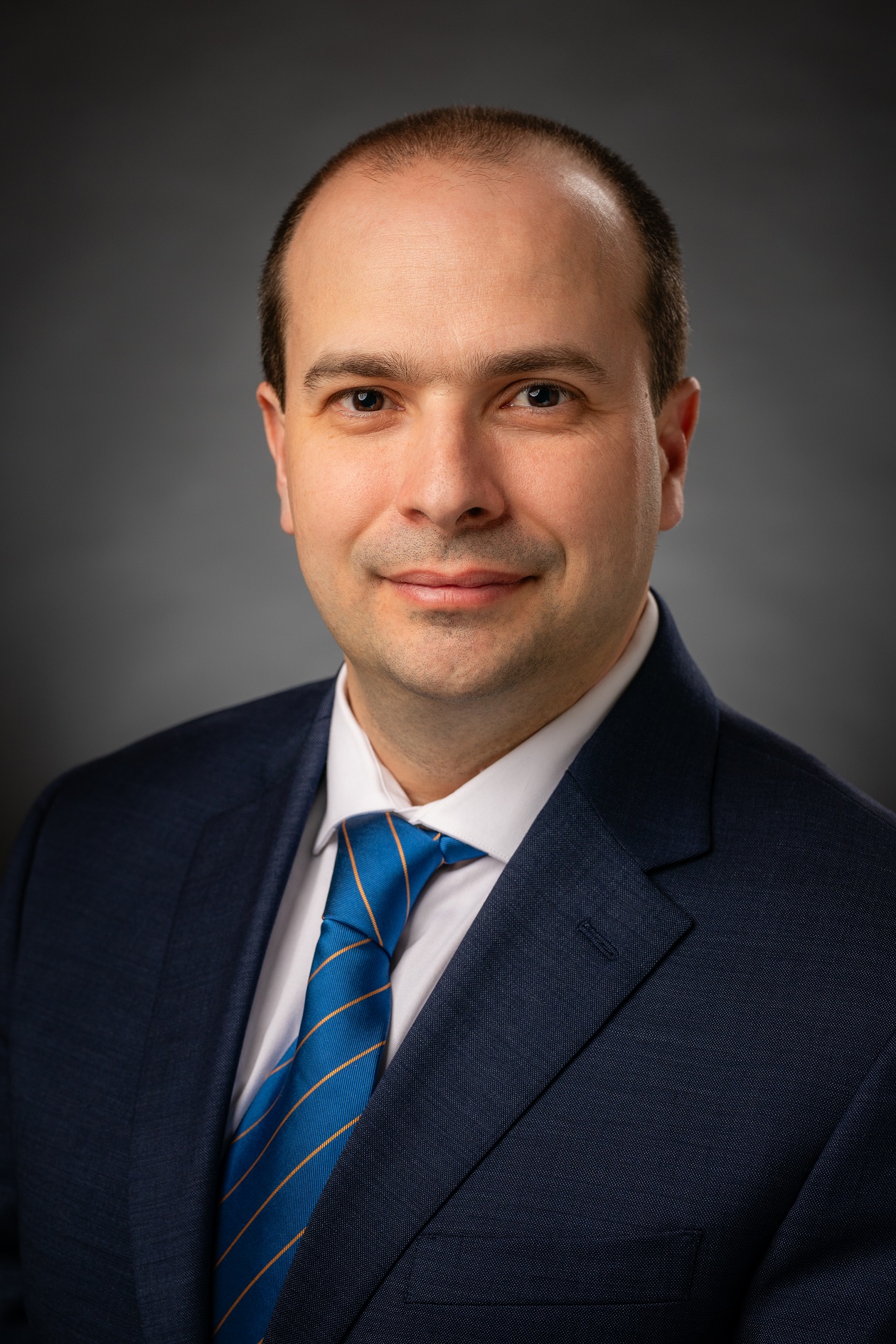}}]{Andrea L'Afflitto}
    (Senior Member, IEEE) is an Associate Professor at the Grado Department of Industrial and Systems Engineering at Virginia Tech, and
    holds affiliate positions with the Department of Mechanical Engineering, the Kevin T. Crofton Department of Aerospace and Ocean Engineering, and the National Security Institute at Virginia Tech.
    Before this appointment, he served as an Assistant Professor at the School of Aerospace and Mechanical Engineering at the University of Oklahoma. 
    Dr. L'Afflitto served as a summer faculty fellow at ARL, AFRL, NAVAIR, and Boeing. 
    His research lies at the intersection of robust adaptive control, native space theory, and hybrid dynamical systems.
    Uncrewed aerial vehicles (UAVs) serve as typical testbeds for numerous of his theoretical results.
    Dr. L'Afflitto is a recipient of the DARPA Young Faculty Award, is an AIAA Associate Fellow, and serves as an Associate Editor-in-Chief for the IEEE Transactions on Aerospace and Electronic Systems and as a member of the IEEE Conference Editorial Board.
\end{IEEEbiography}

\begin{IEEEbiography}[{\includegraphics[width=1in,height=1.25in,clip,keepaspectratio]{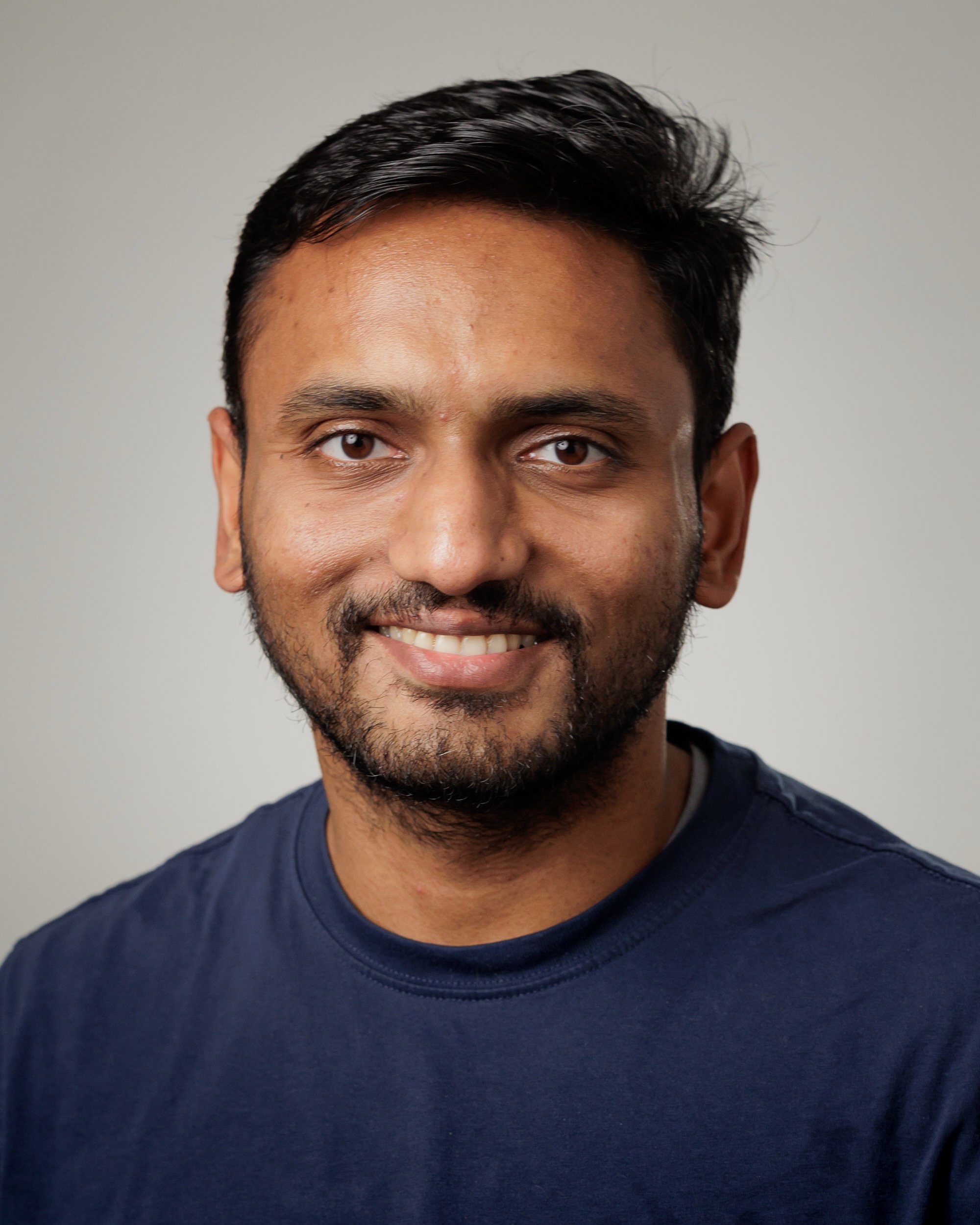}}]{Sai Tej Paruchuri} is a research scientist in the Department of Mechanical Engineering and Mechanics at Lehigh University, Bethlehem, PA, USA. He earned a B.E. in Mechanical Engineering from Thiagarajar College of Engineering, Anna University, Madurai, India, followed by an M.S. in Mathematics and a Ph.D. in Mechanical Engineering from Virginia Tech, Blacksburg, VA, USA. His research interests encompass nonlinear control, adaptive function estimation and control, reproducing kernel Hilbert spaces, data-driven modeling and control, and reinforcement learning. Currently, his research is primarily aimed at addressing control challenges in nuclear fusion and tokamaks.
\end{IEEEbiography}

\begin{IEEEbiography}[{\includegraphics[width=1in,height=1.25in,clip,keepaspectratio]{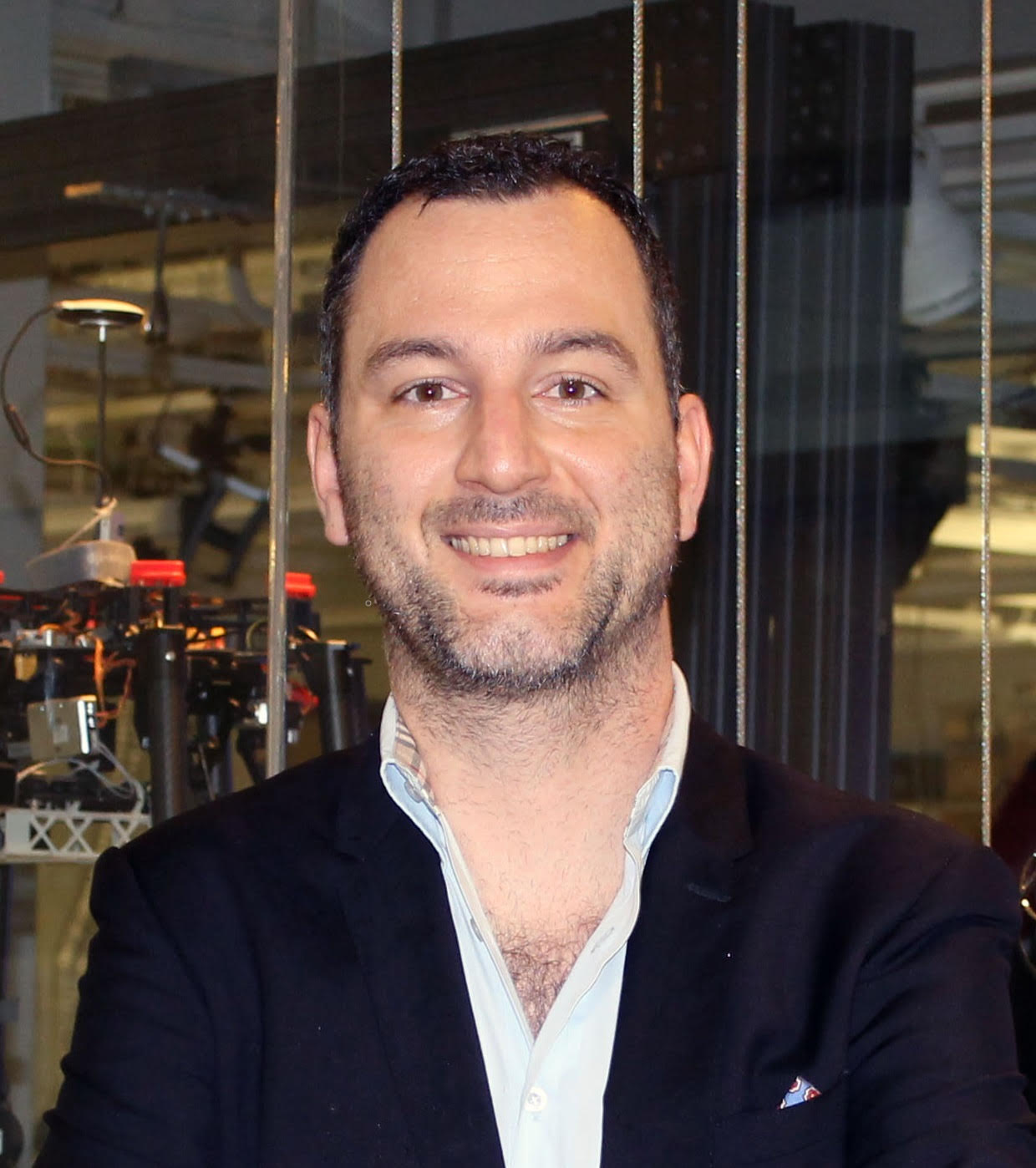}}]{Kyriakos G. Vamvoudakis}
(Senior Member, IEEE) was born in Athens, Greece. He received the Diploma in
Electronic and Computer Engineering from the
Technical University of Crete, Greece in 2006,
and the MSc and PhD degrees in Electrical Engineering at The University of Texas, Arlington in
2008 and 2011, respectively. During the period from 2012 to 2016 he was a project research scientist at the Center for Control, Dynamical Systems and Computation at the University of California, Santa Barbara. He was an Assistant Professor at the Kevin T. Crofton Department of Aerospace and Ocean Engineering at Virginia Tech until 2018.
He currently serves as the Dutton-Ducoffe Endowed Professor at The Daniel Guggenheim School of Aerospace Engineering at Georgia Tech. He holds a secondary appointment in the School of Electrical and Computer Engineering. His expertise is in reinforcement learning, control theory, game theory, cyber-physical security, bounded rationality, and safe/assured autonomy.
He has received numerous prestigious awards, including the 2019 ARO YIP Award, the 2018 NSF CAREER Award, the 2018 DoD Minerva Research Initiative Award, and the 2021 GT Chapter Sigma Xi Young Faculty Award. His work has also been recognized with several best paper nominations and international awards, such as the 2016 International Neural Network Society Young Investigator (INNS) Award.  He is the Editor-in-Chief of Aerospace Science and Technology and currently serves on the IEEE Control Systems Society Conference Editorial Board. Additionally, he is an Associate Editor for several journals, including Automatica, IEEE Transactions on Automatic Control, IEEE Transactions on Neural Networks and Learning Systems, IEEE Transactions on Systems, Man, and Cybernetics: Systems, IEEE Transactions on Artificial Intelligence, Neural Networks, and the Journal of Optimization Theory and Applications. He is also a Senior Guest Editor for the IEEE Open Journal of Control Systems for the special issue on the intersection of machine learning with control. He is a registered Professional Engineer (PE) in Electrical/Computer Engineering, a member of the Technical Chamber of Greece, and an Associate Fellow of AIAA.
\end{IEEEbiography}

\end{document}